\documentclass[10pt]{amsart}
\usepackage{amsmath}
\usepackage{amssymb}
\usepackage{latexsym}
\usepackage{amsbsy}
\usepackage[all]{xy}
\usepackage{amscd}
\usepackage[dvips]{graphicx}

\newtheorem{Prop}{Proposition}[section]
\newtheorem{Thm}[Prop]{Theorem}
\newtheorem{Lemma}[Prop]{Lemma}
\newtheorem{Cor}[Prop]{Corollary}
\newtheorem{Example}[Prop]{Example}

\newtheorem{Definition}[Prop]{Definition}

\def\az{\alpha}      \def\ud{{\underline{d}}}

\def\lz{\lambda}

\def\bbn{{\mathbb N}}  \def\bbz{{\mathbb Z}}  \def\bbq{{\mathbb Q}} \def\bb1{{\mathbb 1}}
\def\mv{{\mathcal V}}
   \def\bbe{{\mathbb E}} 

  \def\bbc{{\mathbb C}}

\def\ra{\rightarrow}
\def\hom{\mbox{Hom}}
\def\rad{\mbox{rad}\,}

\def\ext{\mbox{Ext}\,} 
\def\dim{\mbox{dim}\,}

\def\udim{\mbox{\underline {dim}}}

\def\mod{\mbox{mod}\,}  
    \def\aut{\mbox{Aut}\,}
\def\ker{\mbox{Ker}\,}
\def\cok{\mbox{Coker}\,}

\def\lr#1{\langle #1\rangle}

\def\uq2{U_q(\hat{sl}_2)}

\def\bb{{\bf b}}

\def\nd{{\noindent}}

\def\mc{{\mathcal{C}}}

\def\md{{\mathcal{D}}}
\def\mp{{\mathcal{P}}}

\def\mo{{\mathcal{O}}}

\def\mc{{\mathcal{C}}}

\def\md{{\mathcal{D}}}
\def\mp{{\mathcal{P}}}

\def\mo{{\mathcal{O}}}

\def\mh{{\mathcal{H}}}

\def\ue{{\underline{e}}}

\newcommand{\str}[1]{\langle #1\rangle}

\begin{document}

\title{Green formula in Hall algebras and cluster algebras}

\thanks{ The research was
supported in part by NSF of China (No. 10631010) and by NKBRPC (No. 2006CB805905) \\
2000 Mathematics Subject Classification: 16G20, 16G70. \\ Key words
and phrases: Ringel-Hall algebra, Green's formula,  Cluster
category. }

\author{Jie Xiao and Fan Xu}
\address{Department of Mathematical Sciences\\
Tsinghua University\\
Beijing 10084, P.~R.~China} \email{jxiao@math.tsinghua.edu.cn
(J.Xiao),\  fanxu@mail.tsinghua.edu.cn (F.Xu)}

\maketitle

The objective of the present paper is to give a survey of recent
progress on applications of the approaches of Ringel-Hall type
algebras to quantum groups and cluster algebras via various forms
of Green's formula. In this paper, three forms of Green's formula
are highlighted, (1) the original form of Green's formula
\cite{Green}\cite{RingelGreen}, (2) the degeneration form of
Green's formula \cite{DXX} and (3) the projective form of Green's
formula \cite{XX2007a} i.e. Green formula with a
$\bbc^{*}$-action. The original Green's formula supplies the
comultiplication structure on Ringel-Hall algebras. This
compatibility theorem for multiplication and comultiplication on
Ringel-Hall algebras deals with the symmetric relation between
extensions and flags in the module category of a hereditary
algebra. It provides the quantum group a Hopf algebra structure
and the Drinfeld double in a global way (see \cite{Green},
\cite{Xiao1997} and \cite{Ka}).The second and third section
contribute to these results. The degenerated Green formula can be
viewed as the Green formula over the complex field. We found that
it holds for any algebra given by quiver with relations, not only
hereditary algebras. We write the formula in a geometric version
of the original Green formula, although we know that it
essentially agrees with the restriction functor given by Lusztig
in \cite{Lusztig}. It is applied to provide the geometric
realization of the comultiplication of universal enveloping
algebras. Section 4 is concerned with these results. The
projective version of Green formula has its independent interest.
We give its expression and explain its meaning in Section 6. It is
applied to prove the cluster multiplication theorem which extends
the Caldero-Keller formula \cite{CK2005}. Section 7 is used to
explain the proof in some details. There is a key ingredient
contributing to the above application to cluster categories, i.e.
2-Calabi-Yau property for the cluster categories. In the last
section, we extend the multiplication formula in \cite{GLS2006}
for preprojective algebras to any 2-Calabi-Yau algebras.

Since our main concern in this article is around various forms of
Green formula and their applications, also due to lack of space
and knowledge, we do not include many important topics related to
cluster algebras and Hall algebras. For cluster algebras we refer
to \cite{FZ3} for an excellent survey and we refer to \cite{DX2}
and \cite{Schiffmann} for further study on Ringel-Hall algebras.

\setcounter{tocdepth}{1} \tableofcontents

\section{\bf Representations of quivers and varieties of modules}
\subsection{} 
Let $k$ be a field. A {\it quiver} is a quadruple $Q =
(I,Q_1,s,e)$ where $I$ and $Q_1$ are sets with $I$ non-empty, and
$s,t: Q_1 \to I$ are maps such that $s^{-1}(i)$ and $e^{-1}(i)$
are finite for all $i \in I$. We call $I$ the set of {\it
vertices} and $Q_1$ the set of {\it arrows} of $Q$. For an arrow
$\alpha \in Q_1$ one calls $s(\alpha)$ the starting vertex and
$e(\alpha)$ the end vertex of $\alpha$.

A {\it path} of length $t$ in $Q$ is a sequence $p = \alpha_1
\alpha_2 \cdots \alpha_t$ of arrows such that $s(\alpha_i) =
e(\alpha_{i+1})$ for $1 \leq i \leq t-1$. Set $s(p) = s(\alpha_t)$
and $e(p) = e(\alpha_1)$. Additionally, for each vertex $i \in I$
let $e_i$ be a path of length 0. By $kQ$ we denote the {\it path
algebra} of $Q$, with basis the set of all paths in $Q$ and
product given by concatenation. A ${\it relation}$ for $Q$ is a
linear combination $\sum_{i = 1}^t \lambda_i p_i $ where
$\lambda_i \in k$ and the $p_i$ are paths of length at least two
in $Q$ with $s(p_i) = s(p_j)$ and $t(p_i) = t(p_j)$ for all $1
\leq i,j \leq t$. Thus, we can regard a relation as an element in
$kQ$. An ideal $J$ in $kQ$ is {\it admissible} if it is generated
by a set of relations for $Q$.
\bigskip
\subsection{}
A map $\ud: I \rightarrow \mathbb{N}$ such that $I \setminus
d^{-1}(0)$ is finite is called a {\it dimension vector for} $Q$.
We also write $d_i$ instead of $\ud(i)$, and we often use the
notation $\ud = (d_i)_{i \in I}$. By $\mathbb{N}^{(I)}$ we denote
the semigroup of dimension vectors for $Q$.

A representation $(V,f)$ of the quiver $Q$ over $k$ is a set of
vector spaces $\{V(i)\mid i\in I\}$ together with $k$-linear maps
$f_{\alpha}: V(i)\rightarrow V(j)$ for each arrow $\alpha:
i\rightarrow j.$ For any dimension vector
 $\ud=\sum_i a_i i\in\bbn I,$ we consider the affine space over $k$
$$\bbe_{\ud}(Q)=\bigoplus_{\az\in Q_1}\hom_{k}(k^{a_{s(\az)}},k^{a_{t(\az)}}).$$

For a representation $x = (x_\alpha)_{\alpha \in Q_1} \in
\bbe_{\ud}(Q)$ and a path $p = \alpha_1 \alpha_2 \cdots \alpha_t$
in $Q$ set
$$
x_p = x_{\alpha_1} x_{\alpha_2} \cdots x_{\alpha_t}.
$$
Then $x$ {\it satisfies a relation} $\sum_{i = 1}^t \lambda_i p_i$
if $\sum_{i = 1}^t \lambda_i x_{p_i} = 0$. If $R$ is a set of
relations for $Q$, then let $ \bbe_{\ud}(Q,R) $ be the set of all
representations $x \in \bbe_{\ud}(Q)$ which satisfy all relations
in $R$. This is a closed subvariety of $\bbe_{\ud}(Q)$. Let
$\Lambda$ be the algebra $kQ/J$, where $J$ is the admissible ideal
generated by $R$. We often write $\bbe_{\ud}(A)$ instead of
$\bbe_{\ud}(Q,R).$ Note that every element in
$\bbe_{\ud}(\Lambda)$ can be naturally interpreted as an
$\Lambda$-module structure on the $I$-graded vector space
$k^{\ud}:\oplus_{i\in I}k^{d_i}$ where $\ud=(d_i)_{i\in I}$.
 A {\it dimension vector for} $\Lambda$ is by definition the
same as a dimension vector for $Q$, that is, an element of
$\mathbb{N}^{(I)}$.

 We consider the algebraic group:
$$G_{\ud}:=G_{\ud}(Q)=\prod_{i\in I}GL(a_i,k),$$ which acts on
 $\bbe_{\ud}(Q)$ by $(x_{\az})^{g}_{\az\in Q_1}=(g_{t(\az)}x_{\az}g_{s(\az)}^{-1})_{\az\in Q_1}$ for $g\in G_{\ud}$ and
  $(x_{\az})_{\az\in Q_1}\in\bbe_{\ud}(Q).$
It naturally induces the action of $G_{\ud}$ on $\bbe_{\ud}(Q,R).$
The orbit space is $\bbe_{\ud}(Q,R)/G_{\ud}.$ There is a natural
bijection between the set ${\mathcal M}(\Lambda,\ud)$ of
isomorphism classes of $\bbc$-representations of $\Lambda$ with
dimension vector $\ud$ and the set of orbits of $G_{\ud}$ in
$\bbe_{\ud}(Q,R).$ So we may identify  ${\mathcal M}(\Lambda,\ud)$
with $\bbe_{\ud}(Q,R)/G_{\ud}.$

By abuse of notation, we denote by $\underline{i}$ the dimension
vector mapping $i$ to 1 and $j\neq i$ to $0$. The variety
$\bbe_{\underline{i}}(\Lambda)$ consists just of a single point
corresponding to 1-dimensional $A$-module $S_i$.

An element $x \in \bbe_{\ud}(\Lambda)$ is said to be {\it
nilpotent} if there exists an $N \in \mathbb{Z}^{+}$ such that for
any path $p$ of length greater than $N$ we have $x_p = 0$. By
$\bbe_{\ud}^0(\Lambda)$ we denote the closed subset of nilpotent
elements in $\bbe_{\ud}(\Lambda)$. The nilpotent elements are
exactly the $A$-modules whose composition series contains only
factors isomorphic to the simple modules $S_i$, $i \in I$.
\bigskip
\subsection{}
In the following three subsections, we assume the base field
$k=\bbc.$ The intersection of an open subset and a close subset in
$\bbe_{\ud}(Q,R)$ is called a locally closed subset. A subset in
$\bbe_{\ud}(Q,R)$ is called constructible if and only if it is a
disjoint union of finitely many locally closed subsets. Obviously,
an open set and a close set are both constructible sets. A
function $f$ on $\bbe_{\ud}(Q,R)$ is called constructible if
$\bbe_{\ud}(Q,R)$ can be divided into finitely many constructible
sets satisfies that $f$ is constant on each such constructible
set. Write M(X) for the $\mathbb{Q}$-vector space of constructible
functions on some complex algebraic variety $X$.

Let $\mathcal{O}$ be an above constructible set, $1_{\mathcal{O}}$
is called a characteristic function if $1_{\mathcal{O}}(x)=1$, for
any $x\in \mathcal{O}$; $1_{\mathcal{O}}(x)=0$, for any $x\notin
\mathcal{O}$. It is clear that $1_{\mathcal{O}}$ is the simplest
constructible function and any constructible function is a linear
combination of characteristic functions.  For any constructible
subset $\mathcal{O}$ in $\bbe_{\ud}(Q,R)$, we say $\mathcal{O}$ is
$G_{\ud}$-invariant if $G_{\ud}\cdot\mathcal{O}=\mathcal{O}.$

In the following, the constructible sets and functions will always
be assumed $G_{\ud}$-invariant unless particular statement.
\bigskip
\subsection{} Let $\chi$ denote Euler characteristic in compactly-supported
cohomology. Let $X$ be an algebraic variety and $\mo$ a
constructible subset as the disjoint union of finitely many
locally closed subsets $X_i$ for $i=1,\cdots,m.$ Define
$\chi(\mo)=\sum_{i=1}^m\chi(X_i).$ We note that it is
well-defined. We will use the following properties:
\begin{Prop}[\cite{Riedtmann} and
\cite{Joyce}]\label{Euler} Let $X,Y$ be algebraic varieties over
$\mathbb{C}.$ Then
\begin{enumerate}
    \item  If an algebraic variety $X$ is the disjoint union of
finitely many constructible sets $X_1,\cdots,X_r$, then
$$\chi(X)=\sum_{i=1}^{r}{\chi(X_i)}$$
    \item  If $\varphi:X\longrightarrow Y$ is a morphism
with the property that all fibers have the same Euler
characteristic $\chi$, then $\chi(X)=\chi\cdot \chi(Y).$ In
particular, if $\varphi$ is a locally trivial fibration in the
analytic topology with fibre $F,$ then $\chi(Z)=\chi(F)\cdot
\chi(Y).$
    \item $\chi(\bbc^n)=1$ and $\chi(\mathbb{P}^n)=n+1$ for all $n\geq
    0.$
\end{enumerate}
\end{Prop}
We recall {\it pushforward} functor from the category of algebraic
varieties over $\mathbb{C}$ and the category of
$\mathbb{Q}$-vector spaces (see \cite{Macpherson} and
\cite{Joyce}).

\nd Let $\phi: X\rightarrow Y$ be a morphism of varieties. For
$f\in M(X)$ and $y\in Y,$ define
$$
\phi_{*}(f)(y)=\sum_{c\neq 0}c\chi(f^{-1}(c)\cap \phi^{-1}(y))
$$
\begin{Thm}[\cite{Dimca},\cite{Joyce}]\label{Joyce}
Let $X,Y$ and $Z$ be algebraic varieties over $\mathbb{C},$ $\phi:
X\rightarrow Y$ and $\psi: Y\rightarrow Z$ be morphisms of
varieties, and $f\in M(X).$ Then $\phi_{*}(f)$ is constructible,
$\phi_{*}: M(X)\rightarrow M(Y)$ is a $\mathbb{Q}$-linear map and
$(\psi\circ \phi)_{*}=(\psi)_{*}\circ (\phi)_{*}$ as
$\mathbb{Q}$-linear maps from $M(X)$ to $M(Z).$
\end{Thm}
In fact, Joyce extended this result to algebraic stacks
\cite{Joyce}.
\bigskip
\subsection{} In order to deal with orbit spaces, We also need to consider the
geometric quotients.
\begin{Definition}\label{quotient}
Let $G$ be an algebraic group acting on a variety $X$ and
$\phi:X\rightarrow Y$ be a $G$-invariant morphism, i.e. a morphism
constant on orbits. The pair $(Y,\phi)$ is called a geometric
quotient if $\phi$ is open and for any open subset $U$ of $Y$, the
associated comorphism identifies the ring $\mo_{Y}(U)$ of regular
functions on $U$ with the ring $\mo_{X}(\phi^{-1}(U))^{G}$ of
$G$-invariant regular functions
 on $\phi^{-1}(U)$.
\end{Definition}

The following result due to Rosenlicht \cite{Ro} is essential to
us.

\begin{Lemma}\label{Rosenlicht}
Let $X$ be a $G$-variety, then there exists an  open and dense
$G$-stable subset which has a geometric $G$-quotient.
\end{Lemma}
By this Lemma, we can construct a finite stratification over $X.$
Let $U_1$ be a open and dense $G$-stable subset of $X$ as in Lemma
\ref{Rosenlicht}. Then
$dim_{\mathbb{C}}(X-U_1)<dim_{\mathbb{C}}X.$ We can use the above
lemma again, there exists a dense open $G$-stable subset $U_2$ of
$X-U_1$ which has a geometric $G$-quotient. Inductively, we get
the finite stratification $X=\cup_{i=1}^{l}U_{i}$ which $U_{i}$ is
$G$-invariant locally closed subset and has geometric quotient,
$l\leq dim_{\mathbb{C}}X.$ We denote by $\phi_{U_i}$ the geometric
quotient map on $U_i.$ Define the {\it quasi} Euler-Poincar\'e
characteristic of $X/G$ by
$\chi(X/G):=\sum_{i}\chi(\phi_{U_i}(U_i)).$ If $\{U'_i\}$ is
another choice for the definition of $\chi(X/G)$, then
$\chi(\phi_{U_i}(U_i))=\sum_{j}\chi(\phi_{U_i\cap U'_j}(U_i\cap
U'_j))$ and $\chi(\phi_{U'_j}(U'_j))=\sum_{i}\chi(\phi_{U_i\cap
U'_j}(U_i\cap U'_j)).$ Thus
$\sum_{i}\chi(\phi_{U_i}(U_i))=\sum_{i}\chi(\phi_{U'_i}(U'_i))$
and $\chi(X/G)$ is well-defined (see \cite{XXZ2006}). Similarly,
$\chi(\mo/G):=\sum_i\chi(\phi_{U_i}(\mo\bigcap U_i))$ is
well-defined for any $G$-invariant constructible subset $\mo$ of
$X.$

We also introduce the following notation. Let $f$ be a
constructible function over a variety $X,$ it is natural to define
\begin{equation}\label{integral}
\int_{x\in X}f(x):=\sum_{m\in \bbc}m\chi(f^{-1}(m))
\end{equation}
Comparing with Proposition \ref{Euler}, we also have the following
(see \cite{XXZ2006}).
\begin{Prop}\label{Euler2} Let $X,Y$ be algebraic varieties over
$\mathbb{C}$ under the actions of the algebraic groups $G$ and $H$
respectively.  Then
\begin{enumerate}
    \item  If an algebraic variety $X$ is the disjoint union of
finitely many $G$-invariant constructible sets $X_1,\cdots,X_r$,
then
$$\chi(X/G)=\sum_{i=1}^{r}{\chi(X_i/G)}$$
    \item  If $\varphi:X\longrightarrow Y$ is a morphism induces
    the quotient map $\phi:X/G\rightarrow Y/H$
with the property that all fibers for $\phi$ have the same Euler
characteristic $\chi$, then $\chi(X/G)=\chi\cdot \chi(Y/H).$
   \end{enumerate}
\end{Prop}

Moreover, if there exists an action of algebraic group $G$ on $X$
as Definition \ref{quotient}, and $f$ is $G$-invariant
constructible function over $X,$ We define
\begin{equation}\label{integral2}
\int_{x\in X/G}f(x):=\sum_{m\in \bbc}m\chi(f^{-1}(m)/G)
\end{equation}

In particular, we frequently use the following Corollary.
\begin{Cor}\label{euler3}
Let $X,Y$ be algebraic varieties over $\mathbb{C}$ under the
actions of the algebraic groups $G.$  These actions naturally
induce the action of $G$ on $X\times Y.$ Then
$$\chi(X\times_{G}Y)=\int_{y\in Y/G}\chi(X/G_{y})$$
where $G_y$ is the stable subgroup of $G$ for $y\in Y$ and
$X\times_{G}Y$ is the orbit space of $X\times Y$ under the action
of $G.$
\end{Cor}

\section{\bf Green's formula over finite fields}
In this section, we recall Green's formula  over finite fields
(\cite{Green},\cite{RingelGreen}). Let $k$ be a finite field and
$\Lambda$ a hereditary finitary $k$-algebra. Let $\mathcal{P}$ be
the set of isomorphism classes of finite $\Lambda$-modules. Given
$\alpha\in \mathcal{P},$ let $V_{\alpha}$ be a representative in
$\alpha$ (denoted by $V_{\alpha}\in \alpha$) and
$a_{\alpha}=|\aut_{\Lambda} V_{\alpha}|.$ Given $\xi,\eta$ and
$\lambda$ in $\mathcal{P},$ let $g_{\xi\eta}^{\lambda}$ be the
number of submodules $Y$ of $V_{\lambda}$ such that $Y$ and
$V_{\lambda}/Y$ are isomorphic to $\eta$ and $\xi$, respectively.
The following identity is called Green's formula.
\begin{Thm}\label{Greenformula}
Let $k$ be a finite field and $\Lambda$ a hereditary finitary
$k$-algebra. Let $\xi,\eta,\xi',\eta'\in \mathcal{P}.$ Then
$$
a_{\xi}a_{\eta}a_{\xi'}a_{\eta'}\sum_{\lambda}g_{\xi\eta}^{\lambda}g_{\xi'\eta'}^{\lambda}a_{\lambda}^{-1}
=\sum_{\alpha,\beta,\gamma,\delta}\frac{|\mathrm{Ext}^1(V_{\gamma},V_{\beta})|}{|\mathrm{Hom}(V_{\gamma},V_{\beta})|}g_{\gamma\alpha}^{\xi}
g_{\gamma\delta}^{\xi'}g_{\delta\beta}^{\eta}g_{\alpha\beta}^{\eta'}a_{\alpha}a_{\beta}a_{\delta}a_{\gamma}
$$
\end{Thm}
Given $X\in \xi, Y\in \eta, E\in\lambda,$ let us introduce the
notation
$$h^{\xi\eta}_{\lambda}:=\frac{|\ext^1(X,Y)_{E}|}{|\hom(X,Y)|},$$
where $\ext^1(X,Y)_{E}$ is the subset of $\ext^1(X,Y)$ consisting
of equivalence classes of the exact sequence represented whose
middle term is isomorphic to $E.$ The following identity is called
Riedtmann-Peng formula \cite{Riedtmann}\cite{Peng}.
\begin{Prop}
Let $\alpha,\beta,\lambda\in \mp.$ Then
$$
g_{\alpha\beta}^{\lambda}a_{\alpha}a_{\beta}=h_{\lambda}^{\alpha\beta}
$$
\end{Prop}
Using this proposition, Green's formula can be  rewritten as
(\cite{DXX},\cite{Hubery2005})
\begin{equation}\label{Green-reformulation}
\sum_{\lambda}g_{\xi\eta}^{\lambda}h^{\xi'\eta'}_{\lambda}=\sum_{\alpha,\beta,\gamma,\delta}\frac{|\mathrm{Ext}^1(A,D)||\mathrm{Hom}(M,N)|}{|\mathrm{Hom}(A,D)||\mathrm{Hom}(A,C)||\mathrm{Hom}(B,D)|}g_{\gamma\delta}^{\xi'}g_{\alpha\beta}^{\eta'}h^{\gamma\alpha}_{\xi}h^{\delta\beta}_{\eta}.
\end{equation}
where $X\in \xi, Y\in \eta, M\in\xi', N\in \eta'$ and $A\in
\gamma, C\in \alpha, B\in \delta, D\in \beta, E\in\lambda.$

Let's explain Green's formula by associating some sets to two
sides of the formula. For fixed $kQ$-modules $X,Y,M,N$ with $\udim
X+\udim Y=\udim M+\udim N,$ we fix the $I$-graded $k$-space $E$
such that $ \udim E=\udim X+\udim Y.$ Let $(E,m)$ be a $kQ$-module
structure on $E$ given by the algebraic morphism $m:
\Lambda\rightarrow End_kE.$ Let $Q(E,m)$ be the set of
$(a,b,a',b')$ which satisfies the following diagram with the exact
row and column:
\begin{equation}\label{E:crosses}
\xymatrix{
& & 0 \ar[d] & &\\
& & Y \ar[d]^-{a'} & & \\
0 \ar[r] & N \ar[r]^-{a} & (E,m) \ar[r]^{b} \ar[d]^-{b'} & M \ar[r] & 0\\
& & X \ar[d] & &\\
& & 0 & &}
\end{equation}
Let
$$
Q(X,Y,M,N)=\bigcup_{m:\Lambda\rightarrow End_kE}Q(E,m)
$$
It is clear that
$$
|Q(E,m)|=g_{\xi\eta}^{\lambda}g_{\xi'\eta'}^{\lambda}a_{\xi}a_{\eta}a_{\xi'}a_{\eta'}
$$ where $\lz\in\mathcal{P}$ such that $(E,m)\in\lz,$ or simply write
$m\in\lz.$
$$
|Q(X,Y,M,N)|=|\mathrm{Aut}_{k}E|\cdot\sum_{\lambda}g_{\xi\eta}^{\lambda}g_{\xi'\eta'}^{\lambda}a_{\xi}a_{\eta}a_{\xi'}a_{\eta'}a_{\lambda}^{-1}
$$
It is almost left side of Green's formula.

 Let $\md(X,Y,M,N)^{*}$
be the set of $(B,D,e_1,e_2,e_3,e_4)$ satisfying the following
diagram with exact rows and exact columns:
\begin{equation}\label{E:gurz}
\xymatrix{ & 0 \ar[d] &&&&  0\ar[d] && \\
0 \ar[r] & D \ar[rr]^-{e_1} \ar[dd]^-{u'} && Y\ar[rr]^-{e_2} && B \ar[r] \ar[dd]^-{x}& 0\\
&&&&&&\\
& N \ar[dd]^-{v'} && && M \ar[dd]^-{y} &\\
&&&&&&\\
0\ar[r] &  C \ar[d] \ar[rr]^-{e_3} && X \ar[rr]^-{e_4} && A \ar[r] \ar[d] & 0\\
& 0 && && 0 &}
\end{equation}
where $B,D$ are submodules of $M,N,$ respectively and $A=M/B,
C=N/D.$ The maps $u',v'$ and $x,y$ are naturally induced. We have
$$
|\md(X,Y,M,N)^{*}|=\sum_{\alpha,\beta,\gamma,\delta}g_{\gamma\alpha}^{\xi}g_{\gamma\delta}^{\xi'}g_{\delta\beta}^{\eta}g_{\alpha\beta}^{\eta'}a_{\alpha}a_{\beta}a_{\delta}a_{\gamma}
$$

Fix the above square, let $T=X \times_{A} M =\{(x \oplus m) \in X
\oplus M\;|\; e_4(x)=y(m)\}$ and $S=Y \sqcup_{\small{D}}
\hspace{-.05in}N =Y \oplus N / \{e_1(d)\oplus u'(d)\;|\;d \in
D\}.$ There is a unique map $f:S\rightarrow T$ ( see
\cite{RingelGreen}) such that the natural long sequence
\begin{equation}\label{longexactsequence}
0\rightarrow D\rightarrow S\stackrel{f}\rightarrow T\rightarrow
A\rightarrow 0
\end{equation}
is exact.

We define $(c,d)$ to be the pair of maps satisfying $c$ is
surjective, $d$ is injective and $cd=f.$ The counting of the pair
$(c,d)$ can be made in the following process.  We have the
following commutative diagram:
\begin{equation}\label{3}
\xymatrix{&0\ar[r]& S\ar[r]^{d} \ar[d]& (E,m)
\ar[r]^{d_1}\ar[d]^{c}&
A\ar[r]\ar@{=}[d]&0\\
\varepsilon_0:& 0\ar[r] & Imf \ar[r]& T\ar[r] &A \ar[r]& 0}
\end{equation}
The exact sequence
$$
\xymatrix{0\ar[r]&D\ar[r]& S\ar[r]& Imf \ar[r]&0}
$$
induces the following long exact sequence:
\begin{equation}\label{E:Cross4}
\xymatrix{ 0 \ar[r] & \text{Hom}(A,D) \ar[r] & \text{Hom}(A,S)
\ar[r]& \text{Hom}(A,Imf) \ar[r] &}
\end{equation}
$$\qquad\xymatrix{
 \ar[r] & \text{Ext}^1(A,D) \ar[r] & \text{Ext}^1(A,S) \ar[r]^-{\phi} & \text{Ext}^1(A,Imf) \ar[r] & 0}
$$
We set $\varepsilon_0\in Ext^1(A,Imf)$ corresponding to the
canonical exact sequence
$$
\xymatrix{0\ar[r]&Imf \ar[r]& T \ar[r]& A \ar[r]&0}
$$
and denote $\phi^{-1}(\varepsilon_0)\bigcap
\mathrm{Ext}^1(A,S)_{(E,m)}$ by $\phi^{-1}_{m}(\varepsilon_0).$
 Let $\mathcal{F}(f;m)$ be the set of
$(c,d)$ induced by diagram \eqref{3} with center term $(E,m).$ Let
$$\mathcal{F}(f)=\bigcup_{m:\Lambda\rightarrow
\mathrm{End}_kE}\mathcal{F}(f;m)$$  Then
$$
\mathcal{F}(f;m)=|\phi^{-1}_{m}(\varepsilon_0)|\frac{|\mathrm{Aut}_{\Lambda}(E,m)|}{|\mathrm{Hom}(A,S)|}|\mathrm{Hom}(A,Imf)|
$$
$$
\mathcal{F}(f)=|\mathrm{Aut}_{k}(E)|\frac{|\mathrm{Ext}^1(A,D)|}{|\mathrm{Hom}(A,D)|}
$$

Let $\mo(E,m)$ be the set of $(B,D,e_1,e_2,e_3,e_4,c,d)$
satisfying the following commutative diagram:
\begin{equation}\label{5}
\xymatrix{ & 0 \ar[d] &&&&  0\ar[d] && \\
0 \ar[r] & D \ar[rr]^-{e_1} \ar[dd]^-{u'} && Y\ar[rr]^-{e_2} \ar@{.>}[dl]^-{u_Y}\ar@{.>}[dd]&& B \ar[r] \ar[dd]^-{x}& 0\\
&&S\ar@{.>}[dr]^-{d}&&&&\\
& N \ar@{.>}[ur]_-{u_N}\ar[dd]^-{v'} \ar@{.>}[rr]&& (E,m)\ar@{.>}[rr]\ar@{.>}[dd]\ar@{.>}[dr]^-{c}&& M \ar[dd]^-{y} &\\
&&&&T\ar@{.>}[ur]_-{q_M}\ar@{.>}[dl]^-{q_X}&&\\
0\ar[r] &  C \ar[d] \ar[rr]^-{e_3} && X \ar[rr]^-{e_4} && A \ar[r] \ar[d] & 0\\
& 0 && && 0 &}
\end{equation}
with exact rows and columns, where the maps $q_X,u_Y$ and
$q_M,u_N$ are naturally induced. In fact, the long exact sequence
\eqref{longexactsequence} has the following explicit form:
\begin{equation}\label{longexactsequence2}
\xymatrix{0\ar[r]&
D\ar[r]^{u_Ye_1}&S\ar[r]^{cd}&T\ar[r]^{e_4q_X}&A\ar[r]& 0}
\end{equation}

$$
|\mo(E,m)|=\sum_{\alpha,\beta,\gamma,\delta}|\phi^{-1}_{m}(\varepsilon_0)|\frac{|\mathrm{Aut}_{\Lambda}(E,m)|}{|\mathrm{Hom}(A,S)|}|\mathrm{Hom}(A,Imf)|g_{\gamma\alpha}^{\xi}g_{\gamma\delta}^{\xi'}g_{\delta\beta}^{\eta}g_{\alpha\beta}^{\eta'}a_{\alpha}a_{\beta}a_{\delta}a_{\gamma}
$$
Put $$\mo(X,Y,M,N)=\bigcup_{m:\Lambda\rightarrow
\mathrm{End}_kE}\mo(E,m)$$ Then
$$
|\mo(X,Y,M,N)|=|\mathrm{Aut}_{k}(E)|\sum_{\alpha,\beta,\gamma,\delta}\frac{|\mathrm{Ext}^1(A,D)|}{|\mathrm{Hom}(A,D)|}g_{\gamma\alpha}^{\xi}g_{\gamma\delta}^{\xi'}g_{\delta\beta}^{\eta}g_{\alpha\beta}^{\eta'}a_{\alpha}a_{\beta}a_{\delta}a_{\gamma}
$$
It is almost the right side of Green's formula. Now Green's
formula is equivalent to say
\begin{Prop}\label{Greenbijection}
There exist bijections $\Omega: Q(E,m)\rightarrow \mo(E,m).$
\end{Prop}
By the reformulation \eqref{Green-reformulation} of Green's
formula, Green's formula characterizes the following situation.
Given $\varepsilon\in \ext^1(M,N)_{E}$ and $Y\subseteq E,$ we
count the submodules $M_1\subseteq M, N_1\subseteq N$ induced by
$Y\subseteq E.$ The following proposition deals with the converse
situation when $Y$ is not fixed \cite{Hubery2005}.
\begin{Prop}\label{Huberylemma}
Given $\varepsilon\in \ext^1(M,N)_{E}$ and $M_1\subseteq M,
N_1\subseteq N.$ Set $$V_{M_1N_1}^{\varepsilon}:=\{Y\subseteq
E\mid Y \mbox{ induces }M_1, N_1\}.$$ If it is not a empty set,
then it is isomorphic to $\mathrm{Hom}(M_1,N/N_1).$
\end{Prop}
As Ringel pointed out \cite{RingelGreen}\cite{ZW}, let $\Lambda$
be a finitary $k$-algebra, if Green's formula holds for all finite
$\Lambda$-modules, then the category of finite $\Lambda$-modules
is hereditary. In general, if the category of finite
$\Lambda$-modules is not hereditary, we can extend Green's formula
to non-hereditary case. Given $\Lambda$-modules $A,D,$ let
$\ext^2(A,D)_{S,T}$ denote the subset of $\ext^2_{\Lambda}(A,D)$
consisting of equivalence classes of long exact sequence of the
form
$$
\xymatrix{0\ar[r]&D\ar[r]&S\ar[r]^-{f}&T\ar[r]&A\ar[r]&0}
$$
\begin{Lemma}\label{higherorderextension}
The extension $\varepsilon=0\in Ext^2_{\Lambda}(A,D)$ if and only
$\mathcal{F}(f)\neq \varnothing.$
\end{Lemma}
We write $\varepsilon(\alpha,\beta,\gamma,\delta)=0$ if the
equivalence class in $\mathrm{Ext}^2_{\Lambda}(A,D)$ of the long
exact sequence \eqref{longexactsequence2} vanishes. By Lemma
\ref{higherorderextension}, we have \cite{DXX}
\begin{Thm}
Let $k$ be a finite field and $\Lambda$ a finitary $k$-algebra.
Let $\xi,\eta,\xi',\eta'\in \mathcal{P}.$ Then
$$
a_{\xi}a_{\eta}a_{\xi'}a_{\eta'}\sum_{\lambda}g_{\xi\eta}^{\lambda}g_{\xi'\eta'}^{\lambda}a_{\lambda}^{-1}
=\sum_{\alpha,\beta,\gamma,\delta:
\varepsilon(\alpha,\beta,\gamma,\delta)=0}\frac{|\mathrm{Ext}^1(V_{\gamma},V_{\beta})|}{|\mathrm{Hom}(V_{\gamma},V_{\beta})|}g_{\gamma\alpha}^{\xi}
g_{\gamma\delta}^{\xi'}g_{\delta\beta}^{\eta}g_{\alpha\beta}^{\eta'}a_{\alpha}a_{\beta}a_{\delta}a_{\gamma}.
$$
\end{Thm}
\section{\bf Hall algebras and the realization of quantum groups }
\subsection{}The Ringel-Hall algebra of $\mh(\mathcal{A})$ associated
to an abelian category $\mathcal{A}$ is an associative algebra,
which, as a vector space, has a basis consisting of the
isomorphism classes $[X]$ of objects $X$ in $\mathcal{A},$ and
which has the multiplication $$[X]*[Y]=\sum_{[L]}g_{XY}^{L}[L],$$
where $F_{XY}^{L}$ is the number of subobjects $L'$ of $L$ such
that $L'\cong Y, L/L'\cong X$ and is called the Hall number.
Ringel-Hall algebra has many important variants. We just list here
the following five main types\cite{Joyce-Hall}:
\begin{itemize}
\setlength{\itemsep}{0pt} \setlength{\parsep}{0pt} \item {\bf
Counting subobjects over finite fields}, as in
Ringel~\cite{Ringel}. \item {\bf Perverse sheaves on moduli
spaces} are used by Lusztig~\cite{Lusztig}. \item {\bf Homology of
moduli spaces}, as in Nakajima~\cite{Na}. \item {\bf Constructible
functions on moduli spaces} are used by Lusztig \cite[\S 10.18--\S
10.19]{Lusztig}, Nakajima \cite[\S 10]{Na}, Frenkel, Malkin and
Vybornov \cite{FMV}, Riedtmann \cite{Riedtmann} and others. \item
{\bf Constructible functions on Artin stacks} are used by Joyce
\cite{Joyce-Hall}, Kapranov and Vasserot \cite{KV}.
\end{itemize}
\bigskip
\subsection{} Let $\mathcal A$ be an $k$-linear
abelian category of finite global dimension. The Grothendieck
group $\mathcal G$ of $\mathcal A$ is the quotient of the free
abelian group with generators the isomorphism classes of objects
of $\mathcal A$ by the relations $[X]+[Y]-[L]$ whenever there
exists an exact sequence $0\to Y\to L\to X\to 0$. For any two
objects $M,N\in \mathcal{A},$ we define the multiplicative Euler
form $\langle \,,\, \rangle: \mathcal{G}\times
\mathcal{G}\rightarrow \bbc$ by
$$
\langle M,N\rangle=\sum_{i}(-1)^{i}\mathrm{dim}_{k}\ext^i(M,N)
$$Denote by $h_X$ the image of $[X]$ in $\mathcal G$. We assume
$\mathcal A$ finitary, i.e., each homomorphism group is finite
and, given any $h\in\mathcal G$, there exist only finitely many
isomorphism classes $[X]$ with $h_X=h$. Let $W(X,Y;L)$ denote the
set of all exact sequences $0\to Y\to L\to X\to 0$. The group
$\aut(X,Y):=\aut X\times\aut Y$ acts on $W(X,Y;L)$ via
$$\xymatrix{0\ar[r]&Y\ar[r]^f\ar[d]^\eta &L\ar[r]^g\ar@{=}[d] &X\ar[d]^\xi\ar[r]&0\\
0\ar[r]&Y\ar[r]^{\overline f} &L\ar[r]^{\overline g} &X\ar[r]&0}$$
and we denote the quotient set by $V(X,Y;L)$. Since $f$ is monic
and $g$ epic this action is free, so
$$F_{XY}^L:=\left|V(X,Y;L)\right|=\frac{\left|W(X,Y;L)\right|}{\left|\aut(X,Y)\right|}.$$
By \cite{Riedtmann}, $F_{XY}^{L}=g_{XY}^{L}.$  The Ringel-Hall
algebra $\mathcal H(\mathcal A)$ is the free $\bbz$-module on
generators indexed by the set of isomorphism classes of objects,
writing $u_X$ for $u_{[X]}$,the multiplication is given by
$$u_Xu_Y:=\sum_{[L]}F_{XY}^Lu_L.$$
This sum is finite and that $u_0$ is a unit for the
multiplication.
\begin{Thm}
The Ringel-Hall algebra $\mathcal H(\mathcal A)$ is an
associative, unital algebra.
\end{Thm}
 The associativity of the multiplication is equivalent to say for
 any $X,Y,Z\in \mathcal{A},$ we have
$$\sum_{[L]}F_{XY}^LF_{LZ}^M=\sum_{[L']}F_{XL'}^MF_{YZ}^{L'}$$
\bigskip
\subsection{}Let $k$ be a finite field of order $q_k,$ let
$\Lambda$ be a k-algebra which is hereditary and finitary and
$\mathcal{A}=\mathrm{mod}\Lambda$. Let $R$ be a integral domain
containing the rational field $\bbq$ and an invertible element $v$
such that $v^2=q_k.$ The following we will consider $\mathcal
H(\mathcal A)$ as the free $R$-module with the basis
$\{u_{\alpha}\mid \alpha\in \mp\}.$ Now, we can define the
comultiplication on $\mathcal H(\mathcal A)$ by
$$
\delta: \mathcal H(\mathcal A)\to \mathcal H(\mathcal
A)\otimes_{R} \mathcal H(\mathcal A)
$$
given by
$$
\delta(u_{\lambda})=\sum_{\alpha,\beta}h_{\lambda}^{\alpha\beta}u_{\alpha}\otimes
u_{\beta}.
$$
where $h_{\lambda}^{\alpha\beta}$ is defined in Section 3. We also
consider the tensor product $\mathcal H(\mathcal A)\otimes_{K}
\mathcal H(\mathcal A)$ as an algebra with multiplication $*$ as
follows:
$$
(u_{\alpha}\otimes u_{\beta})*(u_{\gamma}\otimes
u_{\delta})=\frac{|\mathrm{Ext}^1(V_{\alpha},V_{\delta})|}{|\mathrm{Hom}(V_{\alpha},V_{\delta})|}u_{\alpha}u_{u_{\gamma}}\otimes
u_{\beta}u_{\delta}
$$

It is routine to verify that the following fact is an equivalent
version of Green formula in Theorem \ref{Greenformula}. This is
fundamental important to Hall algebras and quantum groups.
\begin{Prop}\label{comultiplication}
The map $\delta$ is a homomorphism between algebras and then
$\mathcal H(\mathcal A)$ is a bialgebra.
\end{Prop}
The twist of the multiplication of $\mathcal H(\mathcal A)$ is
defined by
$$u_Xu_Y:=v^{\langle\udim X, \udim Y\rangle}\sum_{[L]}F_{XY}^Lu_L.$$
With respect to this twist multiplication, we write $\mathcal
H_{*}(\mathcal A)$ instead of $\mathcal H(\mathcal A)$, the
corresponding comultiplication $\delta_{*}$ is defined by
$$
\delta_{*}(u_{\lambda})=\sum_{\alpha,\beta}v^{\langle\udim X,
\udim Y\rangle}h_{\lambda}^{\alpha\beta}u_{\alpha}\otimes
u_{\beta}.
$$

We define a symmetric bilinear form on $\mathcal H(\mathcal A)$
and $\mathcal H_{*}(\mathcal A).$ For $\alpha\in \mp,$ let
$t_{\alpha}=\frac{|V_{\alpha}|}{a_{\alpha}}.$ Defiine
$$
(u_{\alpha},u_{\beta})=\left\{%
\begin{array}{ll}
    t_{\alpha}, & \mbox{ if }\alpha=\beta;\\
    0, & \hbox{otherwise.} \\
\end{array}%
\right.
$$
\begin{Prop}
For any $u_{\alpha},u_{\beta}$ and $u_{\gamma}$ in $\mathcal
H(\mathcal A)$ and $\mathcal H_{*}(\mathcal A),$ we have
$$
(u_{\alpha},u_{\beta}u_{\gamma})=(\delta(u_{\alpha}),u_{\beta}\otimes
u_{\gamma})
$$
and
$$
(u_{\alpha},u_{\beta}*u_{\gamma})=(\delta_{*}(u_{\alpha}),u_{\beta}\otimes
u_{\gamma}).
$$
\end{Prop}
\bigskip
\subsection{} An Euler form is a pair $(\omega,d)$ consisting of a bilinear form $\omega$ on $\bbz I$ with values in $\bbz$ and a function $d: I\to \bbn\cup\{0\}$ such that the following three propeties are satisfied for all $i,j\in I:$
\begin{itemize}
    \item  $\omega(i,j)$ is divisible by $d_i$ and $d_j;$
    \item  $\omega(i,j)\leq 0$ for $i\neq j;$
    \item  $d^{-1}_i\omega(i,i)\leq 1.$
\end{itemize}
The Euler form $(\omega,d)$ is called without short cycles if
$\omega(i,i)=d_i$ and $\omega(i,j)\omega(j,i)=0$ for $i\neq j.$
Given an Euler form without short cycles, its symmetrization will
be a Cartan datum, and any Cartan datum arises in this way
\cite{RingelGreen}. Given a hereditary $k$-algebra, its Euler form
is defined by
$$
\omega(i,j)=\langle S_i,S_j\rangle , d=(d_1,\cdots,d_n)
$$
where $S_i$ are simple $\Lambda$-modules and
$d_i=\mathrm{dim}_{k}\mathrm{End}(S_i).$
\begin{Prop}\label{}
Let $(\omega,d)$ be an Euler form. Let $k$ be a finite field. Then
there exists a hereditary $k$-algebra $\Lambda_k$ with Euler form
$(\omega,d).$
\end{Prop}
Given an Euler form $(\omega,d)$ without short cycles, let
$(I,(-,-))$ be the induced Cartan datum and let $\mathcal{K}$ be a
set of finite fields $k$ such that the set $\{q_k=|k|: k\in
\mathcal{K}\}$ is infinite and $R$ an integral domain containing
$\mathbb{Q},$ and also containing an element $v_k$ such that
$v_k^2=q_k$ for each $k\in \mathcal{K}.$ There exists a hereditary
algebra $\Lambda_k$ with Euler form $(\omega,d).$ Consider the
direct product
$$
\mathcal{H}_{*}(\omega,d):=\prod_{k\in
\mathcal{K}}\mathcal{H}_{*}(\Lambda_k)
$$
By $\mathcal{C}(\omega,d)$ we denote the $\mathbb{Q}$-subalgebra
generated by $t,t^{-1}$ and $u_{i,*}$ whose $k$-components are
$v,v^{-1}$ and $u_{i,k},$ respectively. We define the
$\bbq(t)$-algebra
$$\mathcal{C}_{*}(\omega,d)=\bbq(t)\otimes_{\bbq[t,t^{-1}]}\mathcal{C}(\omega,d).$$
Let $\mathfrak{g}$ be the Kac-Moody algebra of type $(I,(-,-)).$
Furthermore, let $\mathbb{Q}(t)$ be the function field in one
indeterminate $t$ over the field $\bbq$ of rational numbers. The
Drinfeld-Jimbo quantization $U$ of the enveloping algebra
$\mathfrak{U}(\mathfrak{g})$ is by definition the
$\bbq(t)$-algebra generated by $\{E_i,F_i\mid i\in I\}$ and
$\{K_{\mu}\mid \mu\in \bbz I\}$ with quantum Serre relations.
There is a triangular decomposition
$$
U=U^{-}\otimes U^{0}\otimes U^{+}
$$
where $U^{+}$ (resp. $U^{-}$) is the subalgebra generated by $E_i$
(resp. $F_i$), $i\in I$  and $U^{0}$ is the subalgebra generated
by $K_{\mu}, \mu\in \bbz I.$ Depending on Proposition
\ref{comultiplication}, the following theorem gives the
realization of the positive part of a quantum group.
\begin{Thm}
Given an Euler form $(\omega,d)$ without short cycles, let
$(I,(-,-))$ be the induced Cartan datum and $U$ be the
Drinfeld-Jimbo quantum group of type $(I,(-,-)).$ Then the
correspondence $u_{i,*}\mapsto E_i$ induces a Hopf algebra
isomorphism $\mathcal{C}_{*}(\omega,d)\to U^{+}. $
\end{Thm}

\bigskip
\section{\bf Green formula over complex fields with an application to the realization of universal enveloping algebras}
In this section, we assume the base field $k=\bbc.$
\subsection{}Now we introduce the
$\chi$-Hall algebra \cite{Schiffmann}. Assume the base field
$k=\bbc.$ Let $\Lambda=\bbc Q/\langle R\rangle$ be an associative
algebra over $\bbc.$ Let $\mathcal{O}_1,\mathcal{O}_2$ be two
constructible sets with dimension vectors $\underline{d}_1,\ud_2$
invariant under the action of the corresponding group,
respectively,  fix $\Lambda$-module $L$, we define:
$$V(\mathcal{O}_2,\mathcal{O}_1; L)=\{0=X_0\subseteq
X_1\subseteq X_2=L \mid X_i\in \mathrm{mod}\ \Lambda,X_1 \in
\mathcal{O}_1,
 \mbox{  and  } L/X_1\in \mathcal{O}_2\}$$  It
is a constructible subset of Grassmann varieties.  Put
$$
g_{\alpha\beta}^{\lambda}=\chi(\mv(X,Y;L))
$$ for $X\in\alpha, Y\in\beta$  and $L\in\lambda.$ Consider the
$\bbc$-space
$$M_G(\Lambda)=\bigoplus_{\ud\in\bbn I}M_{G_{\ud}}(\Lambda)$$ where
$M_{G_{\ud}}(\Lambda)$ is the $\bbc$-space of $G_{\ud}$-invariant
constructible function on $\bbe_{\ud}(\Lambda).$ We define the
convolution multiplication on $M_G(\Lambda)$ as follows. For any
$f\in M_{G_{\ud}}(\Lambda)$ and $g\in M_{G_{\ud'}}(\Lambda),$
$f\bullet g\in M_{G_{\ud+\ud'}}(\Lambda)$ is given by the formula
$$f\bullet
g(L)=\sum_{c,d\in\bbc}\chi(V(f^{-1}(c),g^{-1}(d);L))cd.$$ We have

\begin{Prop}\label{associativity} The space $M_G(\Lambda)$ under the convolution
multiplication $\bullet$ is an associative $\bbc$-algebra with
unit element.
\end{Prop}

Let $\alpha$ be the image of $X$  in
$\bbe_{\ud_{\alpha}}(Q)/G_{\ud_{\alpha}}.$  We write $X\in
\alpha$, sometimes we also use the notation $\overline{X}$ to
denote the image of $X$ and the notation $V_{\alpha}$ to denote a
representative of $\alpha.$ Instead of $\ud_{\alpha},$ we use
$\underline{\alpha}$ to denote the dimension vector of $\alpha.$
\bigskip
\subsection{}
 For any $L\in \mathrm{mod}\Lambda,$ let
$L=\bigoplus_{i=1}^{r}L_i$ be the decomposition into
indecomposables, then an action of $\bbc^{*}$ on $L$ is defined by
$$
t.(v_1,\cdots,v_r)=(tv_1,\cdots,t^rv_r)
$$
for $t\in \bbc^{*}$ and $v_i\in L_i$ for $i=1,\cdots,r.$

It induces an action of $\bbc^{*}$ on $V(X,Y;L)$ for any
$A$-module $X,Y$ and $L.$ Let $(L_1\subseteq L)\in V(X,Y;L)$ and
$t.L_1$ be the action of $\bbc^{*}$ on $X_1$ by the decomposition
of $L,$ then there is a natural isomorphism $t_{L_1}: L_1\simeq
t.L_1.$ Define $t.(L_1\subseteq L)=(t.L_1\subseteq L).$ We call
this action the first kind of $\bbc^{*}$-action.

Let $D(X,Y)$ be the vector space over $\bbc$ of all tuples
$d=(d(\alpha))_{\alpha\in Q_1}$ such that linear maps
$d(\alpha)\in \mathrm{Hom}_{\bbc}(X_{s(\alpha)},Y_{t(\alpha)})$
and the matrices $
L(d)_{\alpha}=\left(%
\begin{array}{cc}
  Y_{\alpha} & d(\alpha) \\
  0 & X_{\alpha} \\
\end{array}%
\right) $ satisfy the relations $R.$ Define $\pi:
D(X,Y)\rightarrow \mathrm{Ext}^1(X,Y)$ by sending $d$ to the
equivalence class of a short exact sequence
$$
\xymatrix{\varepsilon:\quad 0\ar[r]& Y\ar[rr]^{\left(%
\begin{array}{c}
  1 \\
  0 \\
\end{array}%
\right)}&&L(d)\ar[rr]^{\left(%
\begin{array}{cc}
  0 & 1 \\
\end{array}%
\right)}&&X\ar[r]&0}
$$
where $L(d)$ is the direct sum of $X$ and $Y$ as a vector space.
For any $\alpha\in Q_1,$
$$
L(d)_{\alpha}=\left(%
\begin{array}{cc}
  Y_{\alpha} & d(\alpha) \\
  0 & X_{\alpha} \\
\end{array}%
\right)
$$
Fixed a vector space decomposition $D(X,Y)=\mathrm{Ker}\pi\oplus
E(X,Y),$ then we can identify $\mathrm{Ext}^{1}(X,Y)$ with
$E(X,Y)$ \cite{Riedtmann}\cite{DXX}\cite{GLS2006}. There is a
natural $\bbc^{*}$-action on $E(X,Y)$ by $t.d=(td(\alpha))$ for
any $t\in \bbc^{*}.$ This induces the action of $\bbc^{*}$ on
$\mathrm{Ext}^{1}(X,Y).$ By the isomorphism of $\bbc Q$-modules
between $L(d)$ and $L(t.d),$ we have $t.\varepsilon$is the
following short exact sequence:
$$
\xymatrix{0\ar[r]& Y\ar[rr]^{\left(%
\begin{array}{c}
  t \\
  0 \\
\end{array}%
\right)}&&L(d)\ar[rr]^{\left(%
\begin{array}{cc}
  0 & 1 \\
\end{array}%
\right)}&&X\ar[r]&0}
$$
for any $t\in \bbc^{*}.$ Let $\mathrm{Ext}^{1}(X,Y)_{L}$ be the
subset of $\mathrm{Ext}^{1}(X,Y)$ with the middle term isomorphic
to $L,$ then $\mathrm{Ext}^{1}(X,Y)_{L}$ can be viewed as a
constructible subset of $\mathrm{Ext}^1(X,Y)$ by the
identification between $\mathrm{Ext}^1(X,Y)$ and $E(X,Y).$ Put
$$
h_{\lambda}^{\alpha\beta}=\chi(\mathrm{Ext}^{1}_{A}(X,Y)_{L})
$$
for $X\in \alpha, Y\in \beta$ and $L\in \lambda.$ The following is
known, see a proof in \cite{DXX}. The above $\bbc^{*}$-action on
the extensions induces an action on the middle terms. As a vector
space, $L=X\oplus Y.$ So we can define $t.(x,y)=(x,ty)$ for any
$t\in \bbc^{*}$ and $x\in X, y\in Y.$ Hence, for any $L_1\subseteq
L,$ we obtain a new flag $t.L_1\subseteq L.$ We call this action
the second kind of $\bbc^{*}$-action.
\begin{Lemma}\label{directsum}
For $A,B,X\in mod\Lambda,$
$\chi(\mathrm{Ext}_{\Lambda}^{1}(A,B)_{X})=0$ unless $X\simeq
A\oplus B.$
\end{Lemma}
By considering the $\bbc^{*}$-action on grassmannins and
extensions, we obtain the degeneration form of Green's formula.
\begin{Thm}\label{degenerated}
For fixed $\xi, \eta, \xi', \eta',$ we have
$$
g_{\xi\eta}^{\xi'\oplus
\eta'}=\int_{\alpha,\beta,\delta,\gamma;\alpha\oplus
\gamma=\xi,\beta\oplus \delta=
\eta}g_{\gamma\delta}^{\xi'}g_{\alpha \beta}^{\eta'}
$$
\end{Thm}
This theorem can be induced by consider the morphism
$$
\bigcup_{\alpha,\beta,\gamma,\delta;\alpha\oplus
\gamma=\xi,\beta\oplus \delta= \eta}V(V_{\alpha},V_{\beta};
V_{\eta'})\times V(V_{\gamma},V_{\delta};V_{\xi'})\xrightarrow{i}
V(V_{\xi},V_{\eta};V_{\xi'}\oplus V_{\eta'})
$$
sending $(V^1_{\eta'}\subseteq V_{\eta'},V^1_{\xi'}\subseteq
V_{\xi'})$ to $V^1_{\xi'}\oplus V^1_{\eta'}\subseteq
V_{\xi'}\oplus V_{\eta'}$ in a natural way. We set
$\overline{V}(V_{\xi},V_{\eta};V_{\xi'}\oplus
V_{\eta'}):=V(V_{\xi},V_{\eta};V_{\xi'}\oplus V_{\eta'})\setminus
Im(i).$ This induces the map
$$
\bigcup_{\alpha,\beta,\gamma,\delta; \mathrm{\udim} V_{\beta}+
\mathrm{\udim} V_{\delta}= \ud}V(V_{\alpha},V_{\beta};
V_{\eta'})\times V(V_{\gamma},V_{\delta};V_{\xi'})\xrightarrow{i}
\bigcup_{\mathrm{\udim}
V_{\eta}=\ud}V(V_{\xi},V_{\eta};V_{\xi'}\oplus V_{\eta'})
$$
for fixed $\ud.$

Furthermore, we can consider the inverse morphism:
$$
\bigcup_{\mathrm{\udim}
V_{\eta}=\ud}V(V_{\xi},V_{\eta};V_{\xi'}\oplus
V_{\eta'})\rightarrow
\bigcup_{\alpha,\beta,\gamma,\delta}V(V_{\alpha},V_{\beta};
V_{\eta'})\times V(V_{\gamma},V_{\delta};V_{\xi'})
$$
mapping the submodules of $V_{\xi'\oplus \eta'}$ to the induced
submodules of $V_{\eta'}$ and $V_{\xi'}$ as the bijection in
Proposition \ref{Greenbijection}. This map is a vector bundle by
Proposition \ref{Huberylemma} and is analogous to the restriction
map in \cite{Lusztig}.
\bigskip
\subsection{} A constructible set is called indecomposable if
all points in it correspond to indecomposable $\Lambda$-modules.
Let $\mathcal{O}$ be a constructible set. If it has the form:
$\mo=n_{1}\mathcal{O}_{1}\oplus \cdots \oplus
n_{k}\mathcal{O}_{k}$, where $\mathcal{O}_{i},1\leq i\leq k$ are
indecomposable constructible sets, then $\mo$ is called  to be of
Krull-Schmidt. A constructible set $\mathcal{Q}$ is called to be
of stratified Krull-Schmidt if it has a finite stratification
$\mathcal{Q}=\dot{\bigcup}_{i}\mathcal{Q}_i$ where each
$\mathcal{Q}_i$ is locally closed in $\mathcal{Q}$ and is of
Krull-Schmidt. Define
$R(\Lambda)=\Sigma{\mathbb{Z}1_{\mathcal{Q}}}$, where
$\mathcal{Q}$ are of stratified Krull-Schmidt constructible sets.
Obviously $R(\Lambda)=\Sigma{\mathbb{Z}1_{\mathcal{O}}}$, where
$\mathcal{O}$ are of  Krull-Schmidt constructible sets. We note
that R is not free generated by $1_{\mo}$ for any constructible
subset $\mo.$ By Proposition \ref{associativity},we have
\begin{Thm} The $\bbz$-module $R(\Lambda)$ under the convolution multiplication $\bullet$  is an associative
$\mathbb{Z}$-algebra with unit element.
\end{Thm}

Define the $\mathbb{Z}$-submodule
$L(\Lambda)=\sum_{\mo}{\mathbb{Z}1_{\mathcal{O}}}$ of
$R(\Lambda)$, where $\mathcal{O}$ are indecomposable constructible
sets. Then we have the following result.

\begin{Thm} The $\bbz$-submodule $L(\Lambda)$ is a Lie subalgebra of $R(\Lambda)$ with bracket
$[x,y]=xy-yx$.
\end{Thm}
We consider the following tensor algebra over $\bbz:$
$$T(L)=\bigoplus_{i=0}^{\infty}L^{\otimes i}$$
where $L^{\otimes 0}=\bbz,$ $L^{\otimes i}=L\otimes\cdots\otimes
L$ ( $i$ times ). Then $T(L)$ is an associative $\bbz$-algebra
using the tensor as multiplication. Let $J$ be the two-sided ideal
of $T(L)$ generated by
$$1_{\mo_1}\otimes 1_{\mo_2}-1_{\mo_2}\otimes 1_{\mo_1}-[1_{\mo_1},1_{\mo_2}]$$
where $\mo_1,\ \mo_2$ are any indecomposable constructible sets
and $[-,-]$ is Lie bracket in $L(\Lambda).$ Then
$U(\Lambda)=T(\Lambda)/J$ is the universal enveloping algebra of
$L(\Lambda)$ over $\bbz.$ For a distinguish we write the
multiplication in $U(\Lambda)$ as $\ast$ and the multiplication in
$R(\Lambda)$ as $\bullet.$ We have the canonical homomorphism
$\varphi:U(\Lambda)\rightarrow R(\Lambda)$ satisfying
$\varphi(1_{\mo_1}\ast 1_{\mo_2})=1_{\mo_1}\bullet 1_{\mo_2}.$ Let
$R'(\Lambda)$ be the subalgebra  of $R(\Lambda)$ generated by the
elements
$\lambda_1!\cdots\lambda_n!1_{\lambda_1\mo_1\oplus\cdots\oplus\lambda_n\mo_n},$
where $\lambda_i\in\bbn$ for $i=1,\cdots,n$ and $\mo_i,
i=1,\cdots,n,$ are indecomposable constructible sets and disjoint
to each other.
\begin{Thm} The canonical embedding $\varphi:U(\Lambda)\rightarrow
R(\Lambda)$ induces the isomorphism $\varphi:U(\Lambda)\rightarrow
R'(\Lambda).$ Therefore $U(\Lambda)\otimes_{\bbz}\bbq\cong
R(\Lambda)\otimes_{\bbz}\bbq.$
\end{Thm}
This theorem inspires us that there exists a comultiplication
structure on $R'(\Lambda).$ Indeed, it can be realized by the
degeneration form of Green's formula. Let
$\delta:R'(\Lambda)\longrightarrow
R'(\Lambda)\otimes_{\mathbb{Z}}R'(\Lambda)$ be given by
$\delta(1_{\mathcal{O}_{\lambda}})(A,B)=\chi(\mathrm{Ext}_{\Lambda}^{1}(A,B)_{\mathcal{O}_{\lambda}})$
where $\mathrm{Ext}_{\Lambda}^{1}(A,B)_{\mathcal{O}_{\lambda}}$ is
the set of the equivalence classes of extensions $B$ by $A$ with
middle terms belonging to $\mo_{\lambda}.$
\begin{Thm}
The map $\delta:R'(\Lambda)\longrightarrow R'(\Lambda)\otimes
R'(\Lambda)$ is an algebra homomorphism.
\end{Thm}
\bigskip
\bigskip
\section{\bf Cluster algebras and cluster categories}
\subsection{}Let $\bbc Q$ be a hereditary algebra associated to a
connected quiver without oriented cycles, and let $\md=\md^b(\mod
\bbc Q)$ be the bounded derived category with the shift functor
$M\mapsto M[1]$ for any object $M$ in $\md$. We identify the
category $\mathrm{mod} \bbc Q$ with the full subcategory of $\md$
formed by the complexes whose homology is concentrated in degree
$0$. The Grothendieck group of $\md$ is the same as that of
$\mathrm{mod} \bbc Q$, i.e.,
$\mathcal{G}_Q:=\mathcal{G}(\md)=\mathcal{G}(\mathrm{mod} \bbc Q)$

Let $\tau$ be the Auslander-Reiten translation. It can be
characterized by the Auslander-Reiten formula:
$$
\mathrm{Ext}^{1}_{\md}(M,N)\simeq \mathrm{DHom}_{\md}(N,\tau M)
$$
where $M$, $N$ are any objects in $\md$ and where $D$ is the
functor which takes a vector space to its dual. The AR-translation
$\tau$ is a triangle equivalence and therefore induces an
automorphism of the Grothendieck group of $\md$.

The \emph{cluster category} is defined as the orbit category
$\mc_Q = \md / F$, where $F = \tau^{-1}[1]$ and The category
$\mc_Q$ has the same objects as $\md$, but maps are given by
$$\mathrm{Hom}_{\mc}(X,Y) = \coprod_{i} \mathrm{Hom}_{\md}(F^iX, Y).$$
The category $\mc_Q$ was defined in \cite{CCS} for the A$_n$-case
in \cite{BMRRT} for general case.  Each object $M$ of $\mc_{Q}$
can be uniquely decomposed in the following way:
$$M=M_0\oplus P_M[1],$$
where $M_0$ is the $\bbc Q$-module, and where $P_M[1]$ is the
shift of a projective module.

\begin{Thm}
Let $\bbc Q$ be a hereditary algebra. Let $\mc_Q$ the cluster
category of $\bbc Q$ and $\pi: \md_Q\to \mc_Q$ is the canonical
functor. Then
\begin{itemize}
\item[(a)]{$\mc_Q$ is a Krull-Schmidt category;}
\item[(b)]{$\mc_Q$ is triangulated and $\pi$ is exact;}
\item[(c)]{$\mc_Q$ has AR-triangles and $\pi$ preserves
AR-triangles.}
\end{itemize}
\end{Thm}
The statement (b) is due to Keller \cite{keller2}, while (a) and
(c) are proved in \cite{BMRRT}.

We note that the cluster category is 2-Calabi-Yau which means that
the functor $\mathrm{Ext}^1$ is symmetric in the following sense:
 $$\mathrm{Ext}_{\mc}^1(M,N)\simeq \mathrm{DExt}_{\mc}^1(N,M).$$
\bigskip
\subsection{}
Let $n$ be a positive integer. We fix the {\it ambient field}
${\mathcal{F}}= \bbq(u_1,\ldots,u_m)$ with algebraically
independent generating set $\underline{u}:=(u_1,\ldots,u_n)$. Let
${\mathbf x}$ be a free generating set of $\mathcal{F}$ and let
matrix $B=(b_{ij})_{1\leq i,j\leq n}\in M_n(\bbz)$ is
antisymmetric. Such a pair $({\mathbf x}, B)$ is called {\it a
seed}. The Cartan counterpart of a matrix $B=(b_{ij})$ in
$M_n(\bbz)$ is the matrix $A=(a_{ij})$ with
$$
a_{ij}=\left\{
    \begin{array}{ll}
      2 &\textup{ if } i=j,\\
      -\mid b_{ij}\mid &\textup{ if }i\not=j.
    \end{array}\right.
$$

Let $({\mathbf x},  B)$ be a seed and let $x_j$, $1\leq j\leq n$,
be in ${\mathbf x}$.  We define a new seed as follows. Let $x_j'$
be the element of $\mathcal{F}$ defined by the {\it exchange
relation}:
$$
x_jx_j'=\prod_{b_{ij}>0}x^{b_{ij}}+\prod_{b_{ij}<0}x^{-b_{ij}},
$$
where, by convention, we have $x_i=u_i$ for $i>n$. Set ${\mathbf
x'}={\mathbf x}\cup\{x_j'\}\backslash \{x_j\}$.  Let $ B'$ be the
$n\times n$ matrix given by
$$
b_{ik}'=\begin{cases}-b_{ik}&\hbox{if } i=j \hbox{ or } k=j\\
b_{ik}+\frac{1}{2}( \,|b_{ij}|\, b_{jk}+b_{ij}\, |b_{jk}|\,) &
\hbox{ otherwise.}\cr\end{cases}
$$
Then Fomin and Zelevinsky showed that $({\mathbf x'},  B')$ is
also a seed.  It is called the {\it mutation} of the seed
$({\mathbf x}, B)$ in the direction $x_j$.  Dually, the {\it
mutation} of the seed $({\mathbf x'}, B')$ in the direction $x_j'$
is $({\mathbf x}, B).$ We consider all the seeds obtained by
iterated mutations. The free generating sets occurring in the
seeds are called {\it clusters}, and the variables they contain
are called {\it cluster variables}. By definition, the {\it
cluster algebra} $\mathcal{A}({\mathbf x}, B)$ associated to the
seed $({\mathbf x},B)$ is the $\bbz[u_{1}, \ldots,
u_n]$-subalgebra of $\mathcal{F}$ generated by the set of cluster
variables.

An antisymmetric matrix $B$ defines a quiver $Q=Q_B$ with vertices
corresponding to its rows (or columns) and which has $b_{ij}$
arrows from the vertex $i$ to the vertex $j$ whenever $b_{ij}>0$.
The cluster algebra associated to the seed $({\mathbf x},B)$ will
be also denoted by $\mathcal{A}(Q)$. The following result is
so-called {\it Laurent phenomenon} \cite{FZ2002}.
\begin{Thm}
  Let $B$ be a antisymmetrizable matrix in $M_n(\bbz)$, i.e. there exists a
diagonal matrix in $M_n(\bbn)$ such that $DB$ is antisymmetric.
then the cluster algebra $\mathcal{A}(B)$ associated to $B$ is a
  subalgebra of $\bbq[u_i^{\pm 1}, 1\leq i\leq n]$.
\end{Thm}
A cluster algebra is {\it finite} if it has a finite number of
cluster variables.
\begin{Thm}
  A cluster algebra $\mathcal{A}$ is finite if and only if there exists
  a seed $(\mathbf x, B)$ of $\mathcal{A}$ such that the Cartan counterpart of the
  matrix $B$ is a Cartan matrix of finite type.
\end{Thm}
In this case, there exists a quiver $Q$ of simply laced Dynkin
type associated to the seed $(\mathbf x, B).$

The following proposition is basic for cluster category associated
to a simply laced Dynkin quiver.

\begin{Prop}\cite{BMRRT}\cite{CK2005}\label{clusterbasic}
Let $Q$ be a Dynkin quiver and $M$, $N$ be indecomposable $\bbc
Q$-modules. Let $\mc_Q$ be the corresponding cluster category.
Then
\begin{itemize}
\item[(i)] $\mathrm{Ext}^1_{\mc}(M,N)=\mathrm{Ext}^1_{\bbc
Q}(M,N)\coprod \mathrm{Ext}^1_{\bbc Q}(N,M)$ and at least one of
the two direct factors vanishes. \item[(ii)] any short exact
sequence of $\bbc Q$-modules $\xymatrix{0\ar[r]&M \ar[r]^i & Y
\ar[r]^p & N \ar[r]&0}$ provides a (unique) triangle $\xymatrix{M
\ar[r]^i & Y \ar[r]^p & N \ar[r]& M[1]}$ in $\mc$.
\end{itemize}
\end{Prop}
\bigskip
\subsection{}
An object $M$ of $\mc_Q$ is called {\it exceptional} if it has no
selfextensions, i.e., $\mathrm{Ext}^1_{\mc}(M,M)=0$.  An object
$T$ of $\mc_Q$ is a {\it tilting} object if it is exceptional,
multiplicity free, and has the following maximality property: if
$M$ is an indecomposable object such that $\mathrm{Ext}^1(M,T)=0$,
then $M$ is a direct factor of a direct sum of copies of $T$. Note
that a tilting object can be identified with a maximal set of
indecomposable objects $T_1, \ldots, T_n$ such that
$\mathrm{Ext}^1(T_i,T_j)=0$ for all $i,j$. Define the map
\cite{CCformula}
$$
X_{?}:\mbox{obj}(\mc_Q)\rightarrow \mathbb{Q}(x_1,\cdots,x_n)
$$
satisfying:\\
\nd (i) $X_M$ only depends on the isomorphism class of $M$;\\
\nd (ii) $X_{M\oplus N} = X_M X_N$ for all $M,N$ of $\mc_Q$;\\
\nd (iii) $X_{P_i[1]}=x_i$ for the $ith$ indecomposable projective
$P_i$;\\
\nd (iv) if $M$ is the image in $\mc_Q$ of an indecomposable
$kQ$-module, then
\begin{equation}\label{C-C formula}
X_{M}=\sum_{\ue}\chi(Gr_{\ue}(M))x^{\tau(\ue) -
\underline{\mathrm{dim}}M + \ue}
\end{equation}
where $\tau$ is the Auslander-Reiten translation on the
Grothendieck group $K_0(\md^{b}(Q))$ and, for $v\in \mathbb{Z}^n$,
we put
$$
x^v = \prod_{i=1}^n x_i^{\langle \underline{\mathrm{dim}}S_i,
v\rangle}
$$ and $Gr_{\ue}(M)$ is the $\ue$-Grassmannian of $M,$ i.e. the
variety of submodules of $M$ with dimension vector $\ue.$  This
definition is equivalent to \cite{Hubery2005}
$$
X_{M}=\int_{\alpha,\beta}g_{\alpha\beta}^{M}x^{\underline{\beta}
R+\underline{\alpha} R'-\underline{\mathrm{dim}}M}
$$
where the matrices $R=(r_{ij})$ and $R'=(r'_{ij})$ satisfy
$r_{ij}=\mathrm{dim}_{\bbc}\mathrm{Ext}^1(S_i,S_j)$ and
$r'_{ij}=\mathrm{dim}_{\bbc}\mathrm{Ext}^1(S_j,S_i)$ for $i,j\in
Q_0.$ We note that \cite{Hubery2005}
$$
(\underline{\mathrm{\dim}}{P})R=\underline{\mathrm{\dim}}\rad P
\quad \quad
(\underline{\mathrm{\dim}}{I})R'=\underline{\mathrm{\dim}}I-\underline{\mathrm{\dim}}\mbox{soc}
I
$$
\begin{Prop}
Let $Q$ be a Dynkin quiver and $M$ be a $\bbc Q$ module with
dimension vector $\ud.$ Then
$$
X_M=\frac{1}{x_1^{d_1}x_2^{d_2}\cdots x_n^{d_n}}\sum_{0\leq
\ue\leq \ud}\chi(Gr_{\ue}(M))\prod_{i=1}^n x_i^{\sum_{j\to
i}e_i+\sum_{i\to j}(d_j-e_j)}
$$
\end{Prop}
This formula is called Caldero-Chapoton formula. It induces a
bijection from the set of tilting objects to the set of clusters
of $\mathcal{A}(Q)$.
\begin{Prop}\label{CC-nonproj}
Let $M$ be an indecomposable non-projective $\bbc Q$-module. Then
$$
X_MX_{\tau M}=1+X_E
$$
where $E$ is the middle term of the Auslander-Reiten sequence
ending in $M$.
\end{Prop}
If $M$ is an indecomposable projective $\bbc Q$-module, then
\begin{Prop}\label{CC-proj}
$$
X_{P_i}x_i=1+(X_{\mathrm{rad}P_i})x^{(\mathrm{\udim}S_i)R'},\quad
X_{I_i}x_i=1+(X_{I_i/\mathrm{soc}I_i})x^{(\mathrm{\udim}S_i)R}
$$
where $S_i$ is the $i$th simple module.
\end{Prop}
We give an example to confirm the above properties
\cite{CCformula}.
\begin{Example}
Suppose that $Q$ is the following quiver:
  $$
    \xymatrix{
      1 \ar[r] & 2 & 3 \ar[l]
      }.
  $$
Then, the AR-quiver
  of $\mc_Q$ has the following shape:
  \begin{equation}
    \xymatrix{
       &  P_3[1] \ar[rd] & & P_3 \ar[rd] & & S_1 \ar[dr] & \\
        P_2[1] \ar[ru] \ar[rd] & & S_2=P_2 \ar[ru] \ar[rd] & & I_2 \ar[ru]
       \ar[rd] & &  P_2[1]\\
       &  P_1[1] \ar[ru] & & P_1 \ar[ru] & & S_3 \ar[ur] &.\\
     }
  \end{equation}
  Using formula \eqref{C-C formula}, we can compute explicitly the $X_M$.
  The submodules of $S_2$ are $0$ and $S_2$.
$$
    X_{S_2}=\frac{x_1x_3}{x_2}+\frac{1}{x_2}=\frac{x_1x_3+1}{x_2},
$$
  The submodules of $P_3$ are $0$, $S_2$ and $P_3$ and the submodules
  of $P_1$ are $0$, $S_2$ and $P_1$.
$$
    X_{P_3}=\frac{x_1}{x_2}+\frac{1}{x_2x_3}+\frac{1}{x_3}=\frac{1+x_2+x_1x_3}{x_2x_3},\,
    X_{P_1}=\frac{x_3}{x_2}+\frac{1}{x_2x_1}+\frac{1}{x_1}=\frac{1+x_2+x_1x_3}{x_2x_1},
$$
It is clear that we have
$$
X_{P_3}x_3=1+X_{S_2}, \quad X_{P_1}x_1=1+X_{S_2}
$$
which confirm Proposition \ref{CC-proj}.
  The submodules of $I_2$ are $0$, $P_1$, $P_3$, $S_2$ and $I_2$.
$$
    X_{I_2}=\frac{1}{x_2}+\frac{1}{x_1x_3}+\frac{1}{x_1x_3}+\frac{1}{x_1x_2x_3}+\frac{x_2}{x_1x_3}=
    \frac{1+2x_2+x_2^2+x_1x_3}{x_1x_2x_3},
$$
  The submodules of $S_1$ are $0$ and $S_1$ ; the submodules of $S_3$
  are $0$ and $S_3$.
$$
    X_{S_1}=\frac{1}{x_1}+\frac{x_2}{x_1}=\frac{1+x_2}{x_1},\;
    X_{S_3}=\frac{1}{x_3}+\frac{x_2}{x_3}=\frac{1+x_2}{x_3}.
$$
It is clear that we have
$$
X_{S_1}X_{P_3}=1+X_{I_2}, \quad X_{S_3}X_{P_1}=1+X_{I_2}
$$
which confirm Proposition \ref{CC-nonproj}.
\end{Example}

\section{\bf Green formula under $\bbc^{*}$-action}
\subsection{}
For fixed $\xi,\eta$ and $\xi',\eta'$ with
$\underline{\xi}+\underline{\eta}=\underline{\xi'}+\underline{\eta'}=\underline{\lz}$,
let $V_{\lz}\in\bbe_{\underline{\lz}}$ and $Q(V_{\lz})$ be the set
of $(a,b,a',b')$ which satisfies the following diagram with the
exact row and the exact column:
\begin{equation}\label{E:crosses}
\xymatrix{
& & 0 \ar[d] & &\\
& & V_{\eta} \ar[d]^-{a'} & & \\
0 \ar[r] & V_{\eta'} \ar[r]^-{a} & V_{\lambda} \ar[r]^{b} \ar[d]^-{b'} & V_{\xi'} \ar[r] & 0\\
& & V_{\xi}\ar[d] & &\\
& & 0 & &}
\end{equation}
We let
$$
Q(\xi, \eta, \xi', \eta')=\bigcup_{V_{\lambda}\in
\bbe_{\underline{\lambda}}}Q(V_{\lambda})
$$
We remark $Q(\xi, \eta, \xi', \eta')$ can be viewed as a
constructible subset of module variety
$\bbe_{(\underline{\xi},\underline{\eta},\underline{\xi'},\underline{\eta'},\underline{\lambda})}$
with
$\underline{\xi}+\underline{\eta}=\underline{\xi'}+\underline{\eta'}=\underline{\lambda}$
of the following quiver
\begin{equation}
\xymatrix{
 & 2 \ar[d] & \\
4 \ar[r]& 5 \ar[r] \ar[d] & 3  \\
 & 1 & }
\end{equation}

We have the action of $G_{\underline{\lambda}}$ on $Q(\xi, \eta,
\xi', \eta')$:
$$g.(a,b,a',b')=(ga,bg^{-1},ga',b'g^{-1})$$
The orbit space of $Q(\xi, \eta, \xi', \eta')$ is denoted by
$Q(\xi, \eta, \xi', \eta')^{*}$ and the orbit in $Q(\xi, \eta,
\xi', \eta')^{*}$ is denoted by $(a,b,a',b')^{*}.$ Also we have
the action of $G_{\underline{\lambda}}$ on
$W(V_{\xi'},V_{\eta'};\bbe_{\underline{\lambda}})$ by
$g.(a,b)=(ga,bg^{-1})$ where $$
W(V_{\xi'},V_{\eta'};\bbe_{\underline{\lambda}})=\{(f,g)\mid (f,g)
\mbox{ is an exact sequence with the middle term } L\in
\bbe_{\underline{\lambda}}\}.
$$ This induces the orbit space
$\ext^1(V_{\xi'},V_{\eta'})$ whose orbit is denoted by
$(a,b)^{*}.$ Hence, we have
\begin{equation}\label{2}
\xymatrix{ W(V_{\xi'},V_{\eta'};\bbe_{\underline{\lambda}})\times
W(V_{\xi},V_{\eta};\bbe_{\underline{\lambda}})=Q(\xi, \eta, \xi',
\eta')\ar[d]_-{\phi_1}\ar[r]^-{\phi_2}& Ext^1(V_{\xi'},V_{\eta'})
\\ Q(\xi, \eta,
\xi', \eta')^{*}\ar[ur]_-{\phi}&}
\end{equation}
where $\phi((a,b,a',b')^{*})=(a,b)^{*}$ is well defined.

Let $(a,b,a',b')\in Q(V_{\lambda}),$ then we know \cite{XX2006}
that the stable group of $\phi_1$ at $(a,b,a',b')$ is
$$a'e_1\hom(\cok b'a, \ker ba')e_4b',$$ which is isomorphic to
$\hom(\cok b'a, \ker ba'),$ where the injection $e_1: \ker
ba'\rightarrow V_{\lambda}$ is induced naturally by $a'$ and the
surjection $e_4: V_{\lambda}\rightarrow \cok b'a $ is induced
naturally by $b'$. In the same way, The stable group of $\phi_2$
at $(a,b)$ is $1+a\hom(V_{\xi'},V_{\eta'})b$, which is isomorphic
to $\hom(V_{\xi'},V_{\eta'})$ just for $a$ is injective and $b$ is
surjective. The fibre $\phi_1(\phi_2^{-1}((a,b)^{*}))$ of $\phi$
for $(a,b)^{*}$ has the Euler characteristic
$\chi(W(V_{\xi},V_{\eta};V_{\lambda})).$

Moreover, consider the action of $G_{\underline{\xi}}\times
G_{\underline{\eta}}$ on $Q(\xi,\eta,\xi',\eta')^{*}$ and the
induced orbit space is denoted by $Q(\xi, \eta, \xi',
\eta')^{\wedge}.$  The stable subgroup $G((a,b,a',b')^{*})$ at
$(a,b,a',b')^{*}$ is
$$
\{(g_1,g_2)\in G_{\underline{\xi}}\times G_{\underline{\eta}}\mid
ga'=a'g_2,b'g=g_1b'\quad \mbox{for some }g\in
1+a\hom(V_{\xi'},V_{\eta'})b\}
$$
This determines the group embedding $$G((a,b,a',b')^{*})\ra
(1+a\hom(V_{\xi'},V_{\eta'})b)/(1+ae_1\hom(\cok b'a, \ker
ba')e_4b).$$  The group $G((a,b,a',b')^{*})$ is isomorphic to a
vector space since $ba=0.$ We Know that
$1+a\hom(V_{\xi'},V_{\eta'})b$ is the subgroup of $\aut V_{\lz},$
it acts on $W(V_{\xi},V_{\eta};\bbe_{\underline{\lambda}})$
naturally.
 The orbit space
of $W(V_{\xi},V_{\eta};\bbe_{\underline{\lambda}})$ under the
action of $1+a\mathrm{Hom}(V_{\xi'},V_{\eta'})b$ is denoted by
$\widetilde{W}(V_{\xi},V_{\eta};\bbe_{\underline{\lambda}})$ and a
similar consideration for
$V(V_{\xi},V_{\eta};\bbe_{\underline{\lambda}}).$
 By the above
discussion, we have the following commutative diagram of group
actions:

\begin{equation}\label{stablegroup}
\xymatrix{W(V_{\xi},V_{\eta};\bbe_{\underline{\lambda}})\ar[rrr]^{1+a\hom(V_{\xi'},V_{\eta'})b}
\ar[d]^{G_{\underline{\xi}}\times G_{\underline{\eta}}}&&&
\widetilde{W}(V_{\xi},V_{\eta};\bbe_{\underline{\lambda}})\ar[d]^{G_{\underline{\xi}}\times
G_{\underline{\eta}}}\\
V(V_{\xi},V_{\eta};\bbe_{\underline{\lambda}})\ar[rrr]^{1+a\hom(V_{\xi'},V_{\eta'})b}&&&\widetilde{V}(V_{\xi},V_{\eta};\bbe_{\underline{\lambda}})}
\end{equation}
The stable group for the bottom map is
$$\{g\in 1+a\hom(V_{\xi'},V_{\eta'})b\mid ga'=a'g_2, b'g=g_1b' \mbox{ for some }(g_1,g_2)\in G_{\underline{\xi}}\times G_{\underline{\eta}}\}$$
which is isomorphic to a vector space too, it is denoted by
$V(a,b,a',b').$ We can construct the map from $V(a,b,a',b')$ to
$G((a,b,a',b')^{*})$ sending $g$ to $(g_1,g_2).$ It is
well-defined since $a'$ is injective and $b'$ is surjective. We
have
$$
V(a,b,a',b')/\hom(\cok b'a, \ker ba')\cong G((a,b,a',b')^{*})
$$

We have the following proposition.

\begin{Prop}\label{fibre1} The surjective map
$$\phi^{\wedge}: Q(\xi,
\eta, \xi', \eta')^{\wedge}\rightarrow
\mathrm{Ext}^1(V_{\xi'},V_{\eta'})$$ has fibre isomorphic to
$\widetilde{V}(V_{\xi},V_{\eta};\bbe_{\lambda})$ at $(a,b)^{*}\in
\mathrm{Ext}^1(V_{\xi'},V_{\eta'})_{\lambda}$ where
$\widetilde{V}(V_{\xi},V_{\eta};\bbe_{\lambda})$ satisfies that
there exists a surjective morphism from
$V(V_{\xi},V_{\eta};V_{\lambda})$ to
$\widetilde{V}(V_{\xi},V_{\eta};\bbe_{\lambda})$ with fibre
isomorphic to the vector space
$$\mathrm{Hom}(V_{\xi'},V_{\eta'})/(\mathrm{Hom}(\cok b'a,\ker ba')\times
G((a,b,a',b')^{*})).$$
\end{Prop}
\bigskip
\subsection{} Let $\mo(\xi, \eta, \xi', \eta')$ be the set of
$(V_{\delta},V_{\beta},e_1,e_2,e_3,e_4,c,d)$ satisfying the
following commutative diagram with exact rows and columns:
\begin{equation}\label{guru}
\xymatrix{ & 0 \ar[d] &&&&  0\ar[d] && \\
0 \ar[r] & V_{\beta} \ar[rr]^-{e_1} \ar[dd]^-{u'} && V_{\eta}\ar[rr]^-{e_2} \ar@{.>}[dl]^-{u_{V_{\eta}}}\ar@{.>}[dd]&& V_{\delta} \ar[r] \ar[dd]^-{x}& 0\\
&&S\ar@{.>}[dr]^-{d}&&&&\\
& V_{\eta'} \ar@{.>}[ur]_-{u_{V_{\eta'}}}\ar[dd]^-{v'} \ar@{.>}[rr]&& V_{\lambda}\ar@{.>}[rr]\ar@{.>}[dd]\ar@{.>}[dr]^-{c}&& V_{\xi'} \ar[dd]^-{y} &\\
&&&&T\ar@{.>}[ur]_-{q_{V_{\xi'}}}\ar@{.>}[dl]^-{q_{V_{\xi}}}&&\\
0\ar[r] &  V_{\alpha} \ar[d] \ar[rr]^-{e_3} && V_{\xi} \ar[rr]^-{e_4} && V_{\gamma} \ar[r] \ar[d] & 0\\
& 0 && && 0 &}
\end{equation}
where $V_{\delta},V_{\beta}$ are submodules of
$V_{\xi'},V_{\eta'},$ respectively;
$V_{\gamma}=V_{\xi'}/V_{\delta},$  $
V_{\alpha}=V_{\eta'}/V_{\beta},$ $u', x, v', y$ are the canonical
morphisms, and $V_{\lambda}$ is the center induced by the above
square, $T=V_{\xi} \times_{V_{\gamma}} V_{\xi'} =\{(x \oplus m)
\in V_{\xi} \oplus V_{\xi'}\;|\; e_4(x)=y(m)\}$ and $S=V_{\eta}
\sqcup_{\small{V_{\beta}}} \hspace{-.05in}V_{\eta'} =V_{\eta}
\oplus V_{\eta'} \backslash \{e_1(v_{\beta})\oplus
u'(v_{\beta})\;|\;v_{\beta} \in V_{\beta}\}.$ Then there is unique
map $f:S\rightarrow T$ for the fixed square. We define $(c,d)$ to
be the pair of maps satisfying $c$ is surjective ,$d$ is injective
and $cd=f.$ In particular, its subset with four vertexes
$V_{\gamma},V_{\delta},V_{\alpha},V_{\beta}$ is denoted by
$\mo_{(V_{\gamma},V_{\delta},V_{\alpha},V_{\beta})}.$ There is a
natural action of group $G_{\underline\lz}$ on
$\mo_{(V_{\gamma},V_{\delta},V_{\alpha},V_{\beta})}$ as follows:
$$
g.(V_{\delta},V_{\beta},e_1,e_2,e_3,e_4,c,d)=(V_{\delta},V_{\beta},e_1,e_2,e_3,e_4,cg^{-1},gd)
$$

We denote the orbit spaces of
$\mo_{(V_{\gamma},V_{\delta},V_{\alpha},V_{\beta})}$ and $\mo(\xi,
\eta, \xi', \eta')$  by
$\mo^{*}_{(V_{\gamma},V_{\delta},V_{\alpha},V_{\beta})}$ and
$\mo(\xi, \eta, \xi', \eta')^{*}$ respectively under the actions
of $G_{\underline\lz}.$ The following proposition can be viewed as
the geometrization of Proposition \ref{Greenbijection}.
\begin{Prop}\cite{DXX}There is a homeomorphism
$$
\theta^{*}:Q(\xi, \eta, \xi', \eta')^{*}\rightarrow \mo(\xi, \eta,
\xi', \eta')^{*}
$$
induced by the map between $Q(\xi, \eta, \xi', \eta')$ and
$\mo(\xi, \eta, \xi', \eta').$
\end{Prop}
There is an action of $G_{\underline{\xi}}\times
G_{\underline{\eta}}$ on $\mo(\xi,\eta,\xi',\eta')^{*},$ defined
as follows: for $(g_1,g_2)\in G_{\underline{\xi}}\times
G_{\underline{\eta}},$
$$
(g_1,g_2).(V_{\delta},V_{\beta},e_1,e_2,e_3,e_4,c,d)^{*}=(V_{\delta},V_{\beta},g_2e_1,e_2g_2^{-1},g_1e_3,e_4g_1^{-1},c',d')^{*}
$$
Let us determine the relation between $(c',d')$ and $(c,d).$

Suppose
$(V_{\delta},V_{\beta},g_2e_1,e_2g_2^{-1},g_1e_3,e_4g_1^{-1})$
induces $S',T'$ and the unique map $f': S'\rightarrow T',$ then it
is clear there are isomorphisms:
$$a_1: S\rightarrow S' \quad \mbox{and} \quad a_2: T\rightarrow T'$$
induced by isomorphisms:
$$
\left(%
\begin{array}{cc}
  g_2 & 0\\
  0& id \\
\end{array}%
\right): V_{\eta}\oplus V_{\eta'}\rightarrow V_{\eta}\oplus V_{\eta'} \quad \mbox{and} \quad \left(%
\begin{array}{cc}
  g_1 & 0\\
  0& id \\
\end{array}%
\right): V_{\xi}\oplus V_{\xi'}\rightarrow V_{\xi}\oplus V_{\xi'}
$$

 So $f'=a_2fa_{1}^{-1},$ we have the
following commutative diagram:
\begin{equation}\label{cd}
\xymatrix{S'\ar[r]^{d'}& V_{\lambda}\ar[r]^{d'_1}& V_{\gamma}\ar@{=}[d]\\
S\ar[r]^{d}\ar[u]^{a_1}\ar[d]& V_{\lambda}\ar[r]^{d_1}\ar[d]^{c}\ar[u]^{g}& V_{\gamma}\ar@{=}[d]\\
Imf\ar[r]\ar[d]^{a_2}& T\ar[r]\ar[d]^{a_2}& V_{\gamma}\ar@{=}[d]\\
Imf'\ar[r]& T'\ar[r]& V_{\gamma}}
\end{equation}
Hence, $c'=a_2cg^{-1}$ and $d'=gda_{1}^{-1}.$ In particular,
$c=c'$ and $d=d'$ if and only if $g_1=id_{V_{\xi}}$ and
$g_2=id_{V_{\eta}}.$ This shows the action of
$G_{\underline{\xi}}\times G_{\underline{\eta}}$ is free.

Its orbit space is denoted by $\mo(\xi,\eta,\xi',\eta')^{\wedge}.$
The above homeomorphism induces the Proposition
\begin{Prop}\label{fibre2}
There exists a homeomorphism under quotient topology
$$
\theta^{\wedge}:Q(\xi,\eta,\xi',\eta')^{\wedge}\rightarrow
\mo(\xi,\eta,\xi',\eta')^{\wedge}.
$$
\end{Prop}

\bigskip
\subsection{}Let $\md(\xi,\eta,\xi',\eta')^{*}$ be the set of
$(V_{\delta},V_{\beta},e_1,e_2,e_3,e_4)$ satisfying the diagram
\eqref{guru}, in particular, its subset with four vertexes
$V_{\gamma},V_{\delta},V_{\alpha},V_{\beta}$ is denoted by
$D^{*}_{(V_{\gamma},V_{\delta},V_{\alpha},V_{\beta})}.$ Then we
have a projection:
$$
\varphi^{*}:\mo(\xi,\eta,\xi',\eta')^{*}\rightarrow
\md(\xi,\eta,\xi',\eta')^{*}.
$$
The fibre of this morphism is a affine space with dimension
$\mathrm{dim}_{\bbc}\mathrm{Ext}^1(V_{\gamma},V_{\beta})$ for any
element in
$\md_{(V_{\gamma},V_{\delta},V_{\alpha},V_{\beta})}^{*}.$

There is also an action of group $G_{\underline{\xi}}\times
G_{\underline{\eta}}$ on $\md(\xi,\eta,\xi',\eta')^{*}$ with the
stable group isomorphic to the vector space
$\hom(V_{\gamma},V_{\alpha})\times \hom(V_{\delta},V_{\beta}).$
The orbit space is denoted by $\md(\xi,\eta,\xi',\eta')^{\wedge}.$
The projection naturally induces the projection:
$$
\varphi^{\wedge}: \mo(\xi,\eta,\xi',\eta')^{\wedge}\rightarrow
\md(\xi,\eta,\xi',\eta')^{\wedge}
$$
Its fibre for $(V_{\delta},V_{\beta},e_1,e_2,e_3,e_4)^{\wedge}$ is
isomorphism to the quotient space of
$$(\varphi^{*})^{-1}(V_{\delta},V_{\beta},e_1,e_2,e_3,e_4)$$ under the action of
$\hom(V_{\gamma},V_{\alpha})\times \hom(V_{\delta},V_{\beta}).$
The corresponding stable subgroup is
$$
\{(g_1,g_2)\in 1+e_3\mathrm{Hom}(V_{\gamma},V_{\alpha})e_4\times
1+e_1\mathrm{Hom}(V_{\delta},V_{\beta})e_2\mid
$$$$ga'=a'g_2,b'g=g_1b'\quad \mbox{for some }g\in 1+a\mathrm{Hom}(V_{\xi'},V_{\eta'})b\}
$$
where $a,b,a',b'$ are induced by diagram \eqref{guru}. It is
isomorphic to the vector space $G((a,b,a',b')^{\wedge}).$
Therefore, we have
\begin{Prop}\label{fibre3}
There exists a projection $$ \varphi^{\wedge}:
\mo(\xi,\eta,\xi',\eta')^{\wedge}\rightarrow
\md(\xi,\eta,\xi',\eta')^{\wedge}
$$
with the fibre isomorphic to affine space of dimension
$$
\mathrm{dim}_{\bbc}(\mathrm{Ext}^1(V_{\gamma},V_{\beta})\times
G((a,b,a',b')^{*})/\mathrm{Hom}(V_{\gamma},V_{\alpha})\times
\mathrm{Hom}(V_{\delta},V_{\beta})).
$$
\end{Prop}
Now we have the following diagram of morphisms:
\begin{equation}\label{picture}
\xymatrix{\mathrm{Ext}^1(V_{\xi'},V_{\eta'})&
Q(\xi,\eta,\xi',\eta')^{\wedge}\ar[l]_-{\phi^{\wedge}}\ar[r]^-{\theta^{\wedge}}&\mo(\xi,\eta,\xi',\eta')^{\wedge}\ar[r]^-{\varphi^{\wedge}}&
\md(\xi,\eta,\xi',\eta')^{\wedge}}
\end{equation}
\bigskip
\subsection{} In this subsection, we consider the action of $\mathbb{C}^{*}$ on each term in
\eqref{picture}.

\nd (1) For $t\in \bbc^{*}$ and $\varepsilon=(a,b)^{*}\in
\mathrm{Ext}^1(V_{\xi'},V_{\eta'}),$ set
$t.\varepsilon=(t^{-1}a,b)^{*}.$

\nd (2) For $t\in \bbc^{*}$ and $(a,b,a',b')^{\wedge}\in
Q(\xi,\eta,\xi',\eta')^{\wedge},$ we note that
$Q(\xi,\eta,\xi',\eta')^{\wedge}$ is just
$$
\{(\varepsilon, L_1\subseteq L)\mid\varepsilon\in
\mathrm{Ext}^1(V_{\xi'},V_{\eta'})_{L}, L_1\cong V_{\eta},
L/L_1\cong V_{\xi}\}.
$$
Set
$$t.(a,b,a',b')^{\wedge}=(t^{-1}a,b,t_{a'(V_{\eta})}a',b't^{-1}_{b'})^{\wedge}$$
where $$t_{a'(V_{\eta})}: a'(V_{\eta})\rightarrow t.a'(V_{\eta})$$
and $$t_{b'}: V_{\lambda}/a'(V_{\eta})\rightarrow
V_{\lambda}/t.a'(V_{\eta})$$ are under the second kind of
$\bbc^{*}$-action in Section 5.2.

\nd (3) For $t\in \bbc^{*}$ and
$(V_{\delta},V_{\beta},e_1,e_2,e_3,e_4,c,d)^{\wedge}\in
\mo(\xi,\eta,\xi',\eta')^{\wedge},$ set
$$t.(V_{\delta},V_{\beta},e_1,e_2,e_3,e_4,c,d)^{\wedge}=(t.V_{\delta},t.V_{\beta},t^{-1}e_1t^{-1}_{V_{\beta}},t_{V_{\delta}}
e_2,t^{-1}e_3t^{-1}_{V_{\alpha}},t_{V_{\gamma}}e_4,E,t.c,t.d)^{\wedge}.$$

\nd (4) For $t\in \bbc^{*}$ and
$(V_{\delta},V_{\beta},e_1,e_2,e_3,e_4)^{\wedge}\in
\md(\xi,\eta,\xi',\eta')^{\wedge},$ set
$$t.(V_{\delta},V_{\beta},e_1,e_2,e_3,e_4)^{\wedge}=
(t.V_{\delta},t.V_{\beta},t^{-1}e_1t^{-1}_{V_{\beta}},t_{V_{\delta}}e_2,t^{-1}e_3t^{-1}_{V_{\alpha}},t_{V_{\gamma}}
e_4)^{\wedge}$$ where $t_{V_{\alpha}}=a^{-1}t_{b'}a,
t_{V_{\beta}}=a^{-1}t_{a'(V_{\eta})}a,
t_{V_{\gamma}}=bt_{b'}b^{-1},
t_{V_{\delta}}=bt_{a'(V_{\eta})}b^{-1}$ and $t.V_{\delta},
t.V_{\beta}$ are under the first kind of $\bbc^{*}$-action in
Section 5.2.

By considering the Euler characteristic of the orbit space of
every term in diagram \ref{picture} under $\bbc^{*}$-action, we
have the following theorem, which can be viewed as a geometric
version of Green's formula under the $\bbc^*$-action.
\begin{Thm}\label{projgreen}
For fixed $\xi, \eta, \xi', \eta',$ we have
\begin{eqnarray*}
 &&\int_{\lambda\neq \xi'\oplus
\eta'}\chi(\mathbb{P}Ext^{1}(V_{\xi'},V_{\eta'})_{\lambda})g_{\xi\eta}^{\lambda}=\\
   &&\int_{\alpha,\beta,\delta,\gamma,\alpha\oplus \gamma\neq \xi
\mbox{ or }\beta\oplus \delta\neq
\eta}\chi(\mathbb{P}(Ext^{1}(V_{\gamma},V_{\alpha})_{\xi}\times
Ext^1(V_{\delta},V_{\beta})_{\eta}))g_{\gamma\delta}^{\xi'}g_{\alpha\beta}^{\eta'}\\
   &+& \int_{\alpha,\beta,\delta,\gamma,\alpha\oplus \gamma=\xi,
\beta\oplus\delta=
\eta}[d(\xi',\eta')-d(\gamma,\alpha)-d(\delta,\beta)-\lr{\gamma,\beta}]g_{\gamma\delta}^{\xi'}g_{\alpha\beta}^{\eta'}\\
   &-& \int_{\alpha,\beta,\delta,\gamma,\alpha\oplus \gamma=\xi,
\beta\oplus\delta=
\eta}\chi(\mathbb{P}\overline{V}(V_{\xi},V_{\eta};V_{\xi'}\oplus
   V_{\eta'}))
\end{eqnarray*}
where $\mathbb{P}\overline{V}(V_{\xi},V_{\eta};V_{\xi'}\oplus
   V_{\eta'})$ is the orbit space of $\overline{V}(V_{\xi},V_{\eta};V_{\xi'}\oplus
   V_{\eta'})$ defined in Section 5.2.
   \end{Thm}

\section{\bf Cluster multiplication formula}
\subsection{} Let $Q$ be a simply-laced Dynkin quiver. Let $\mc_Q$ and $\mathcal{A}(Q)$ be the corresponding cluster
category and cluster algebra respectively. As Section 6 showed,
$\mathcal{A}(Q)$ is a cluster categories of finite type. A natural
question is how to realize the cluster algebra $\mathcal{A}(Q)$ by
the cluster category $\mc_Q.$ In \cite{CK2005}, the authors
extended the Caldero-Chapoton formula to give a cluster
multiplication formula to answer this question. For a variety $X$,
we define $\chi_c(X)$ to be the Euler-Poincar\'e characteristic of
the \'etale cohomology with proper support of $X.$
\begin{Thm}\label{CKtheorem}
For any objects $M$, $N$ of $\mc$, we have
\begin{itemize}
\item[(i)] If $\mathrm{Ext}^1(M,N)=0$, then $X_MX_N=X_{M\oplus
N}$, \item[(ii)] If $\mathrm{Ext}^1(M,N)\neq 0$, then
$$
\chi_c(\mathbb{P}\mathrm{Ext}^1(M,N)) X_M X_N = \sum_Y
(\chi_c(\mathbb{P}\mathrm{Ext}^1(M,N)_Y) +
\chi_c(\mathbb{P}\mathrm{Ext}^1(N,M)_Y)) X_Y ,
$$
where $Y$ runs through the isoclasses of $\mc$.
\end{itemize}
\end{Thm}
The formula in this theorem is called Caldero-Keller formula. The
proof of this theorem depends on the Caldero-Chapoton formula and
the fact that the cluster category is 2-Calabi-Yau described as
follows under the context of module category. For $M_1\subseteq M$
and $N_1\subseteq N,$ we consider the map
$$
\beta': \ext^{1}(M,N)\oplus \ext^{1}_{\Lambda}(M_1,N_1)\rightarrow
\ext^{1}_{\Lambda}(M_1,N)
$$
sending $(\varepsilon,\varepsilon')$ to
$\varepsilon_{M_1}-\varepsilon'_{N}$ where $\varepsilon_{M_1}$ and
$\varepsilon'_{N}$ are induced by including $M_1\subseteq M$ and
$N_1\subseteq N,$ respectively as follows:
$$
\xymatrix{
\varepsilon_{M_1}:&0\ar[r]&N\ar[r]\ar@{=}[d]&L_1\ar[r]\ar[d]&M_1\ar[r]\ar[d]&0\\
\varepsilon: & 0\ar[r]&N\ar[r]&L\ar[r]^-{\pi}&M\ar[r]&0 }
$$
where $L_1$ is the pullback, and
$$
\xymatrix{
\varepsilon':&0\ar[r]&N_1\ar[r]\ar[d]&L'\ar[r]\ar[d]&M_1\ar[r]\ar@{=}[d]&0\\
\varepsilon'_N: & 0\ar[r]&N\ar[r]&L'_1\ar[r]&M_1\ar[r]&0 }
$$
where $L'_1$ is the pushout. It is clear that
$\varepsilon,\varepsilon'$ and $M_1,N_1$ induce the inclusions
$L_1\subseteq L$ and $L'\subseteq L'_1.$ Set Let
$$
p_0: \ext^{1}_{\Lambda}(M,N)\oplus
\ext_{\Lambda}^{1}(M_1,N_1)\rightarrow \ext^{1}_{\Lambda}(M,N)
$$
Using 2-Calabi-Yau property (Auslander-Reiten formula)
$\ext^{1}(M,N)\simeq D\hom(N,\tau M),$ we can consider the dual of
$\beta'$
$$
\beta: \hom(N,\tau M_1)\rightarrow \hom(N,\tau M)\oplus
\hom(N_1,\tau M_1)
$$
Using bilinear form and orthogonality, we know
\begin{Prop}
$$
(p_0(Ker\beta'))^{\perp}= Im\beta\bigcap \hom(N,\tau M).
$$
\end{Prop}
We give an example to illustrate Theorem \ref{CKtheorem}
\cite{CK2005}.
\begin{Example}
  Assume $Q$ is the quiver of type A$_2$:
 $$
    \xymatrix{1 & 2\ar[l]}.
$$
 Set $M=S_2\oplus S_2$, $N=S_1\oplus S_1$. The middle term $Y$ is either $S_1\oplus P_2\oplus S_2$ or $P_2\oplus P_2$ if $Y$ is an object such that $\mathbb{P}\mathrm{Ext}^1(M,N)_Y$ is not empty.
  The cardinality of $\mathbb{P}\mathrm{Ext}^1(N,M)_Y$ on $\mathbb{F}q$ is respectively $q^2+2q+1$ and $q(q^2-1)$. In a dual way, the middle term $Y$ is either $S_1\oplus S_2$ or $0$ if  $\mathbb{P}\mathrm{Ext}^1(N,M)_Y$ is not empty.
 The cardinality of $\mathbb{P}\mathrm{Ext}^1(N,M)_Y$ on $\mathbb{F}q$ is respectively $q^2+2q+1$ and $q(q^2-1)$.
 The cluster multiplication theorem gives:
 $$X_NX_M=X_{S_1\oplus P_2\oplus S_2}+X_{S_1\oplus S_2}.$$
\end{Example}
\bigskip
\subsection{} Let $\Lambda=\bbc Q$ be a hereditary algebra associated to a connected quiver without oriented cycles. The following theorem generalizes Caldero-Keller's cluster multiplication theorem for cluster categories of finite type to the following theorem for cluster categories of any
type. Define $$\hom(L_1,L_2)_{Y[1]\oplus
X}=\{g\in\hom(L_1,L_2)|\ker g\simeq Y, \cok g\simeq X\}$$
\begin{Thm}\label{clustertheorem} (1) For any
$\Lambda$-modules $V_{\xi'}$, $V_{\eta'}$ we have
$$\hspace{0cm}d^1(\xi', \eta')X_{V_{\xi'}} X_{V_{\eta'}} =\int_{\lambda\neq
\xi'\oplus \eta'} \chi(\mathbb{P}\mathrm{Ext}^1(V_{\xi'},
V_{\eta'})_{V_{\lambda}})
X_{V_{\lambda}}$$$$+\int_{\gamma,\beta,\iota}\chi(\mathbb{P}\mathrm{Hom}(V_{\eta'},\tau
V_{\xi'})_{V_{\beta}[1]\oplus \tau V_{\gamma}'\oplus
I_0})X_{V_{\gamma}}X_{V_{\beta}}x^{\underline{\mathrm{dim}}
\mathrm{soc}I_0}
$$
where  $I_0\in \iota$ is injective and
$V_{\gamma}=V_{\gamma}'\oplus P_0,$ $P_0$ is the projective direct
summand of $V_{\xi'}.$

(2) For any $\Lambda$-module $V_{\xi'}$ and  $P\in \rho$ is
projective Then
 $$d(\rho, \xi')X_{V_{\xi'}}x^{\underline{\mathrm{dim}}P/radP}=\int_{\delta,\iota'}\chi(\mathbb{P}\mathrm{Hom}(V_{\xi'},I)_{V_{\delta}[1]\oplus
I'})X_{V_{\delta}}x^{\underline{\mathrm{dim}}\mathrm{soc}I'}$$$$+\int_{\gamma,\rho'}\chi(\mathbb{P}\mathrm{Hom}(P,V_{\xi'})_{P'[1]\oplus
V_{\gamma}})X_{V_{\gamma}}x^{\underline{\mathrm{dim}}P'/radP'}
 $$
where $I=DHom(P,\Lambda),$ and $I'\in \iota'$ injective, $P'\in
\rho'$ projective.
\end{Thm}
\nd {\bf Remark.} The formula in this theorem varies from the
reformulation in \cite[Theorem 12]{Hubery2005}.  It extends
Theorem \ref{CKtheorem}. Indeed, if $\Lambda=\bbc Q$ for a simply
laced Dynkin quiver, then by Proposition \ref{clusterbasic}, we
can assume $\mathrm{Ext}^1_{\Lambda}(N,M)=0$ and
$\mathrm{Ext}_{\mc_Q}(M,N)=\mathrm{Ext}^1_{\Lambda}(M,N).$
Furthermore, we have \cite[Lemma 13]{Hubery2005}
$$
\mathrm{Ext}_{\mc_Q}(M,N)_E\simeq \mathrm{Ext}^1_{\Lambda}(M,N)_E
$$
and
$$\mathrm{Ext}_{\mc_Q}(N,M)_{K\oplus C[-1]}\simeq \mathrm{Hom}(N,\tau
M)_{K[1]\oplus C}
$$

The proof of this theorem depends on Theorem \ref{projgreen} and
the projectivization of ``Higher order'' associativity described
as follows. For any $\Lambda$-modules $X,Y,L_1$ and $L_2,$ we
define
$$
\hspace{-9cm}W(X,Y;L_1,L_2):=$$$$\{(f,g,h)\mid
\xymatrix{0\ar[r]&Y\ar[r]^{f}&L_1\ar[r]^{g}&L_2\ar[r]^{h}&
X\ar[r]& 0}~~\mbox{is an exact sequence}\} .$$ Under the actions
of $G_{\underline{\alpha}}\times G_{\underline{\beta}},$ where
$\underline{\alpha}=\udim X$ and $\underline{\beta}=\udim Y,$ the
orbit space is denoted by $V(X,Y;L_1,L_2).$ In fact,
$$
V(X,Y;L_1,L_2)=\{g:L_1\rightarrow L_2\mid \mathrm{Ker}g\cong Y
\mbox{ and } \mathrm{Coker}g\cong X\}
$$
Put
$$
h_{XY}^{L_1L_2}=\chi(V(X,Y;L_1L_2))
$$
We have the following ``higher order'' associativity.
\begin{Thm}\label{associativity}
For fixed $\Lambda$-modules $X,Y_i,L_i$ for $i=1,2,$ we have
$$
\int_{\overline{Y}}g_{Y_2Y_1}^{Y}h_{XY}^{L_1L_2}=\int_{\overline{L'}_1}g_{L'_1Y_1}^{L_1}h_{XY_2}^{L'_1L_2}
$$
Dually, for fixed $\Lambda$-modules $X_i,Y,L_i$ for $i=1,2,$ we
have
$$
\int_{\overline{X}}g_{X_2X_1}^{X}h_{XY}^{L_1L_2}=\int_{\overline{L'}_2}g_{X_2L'_2}^{L_2}h_{X_1Y}^{L_1L'_2}
$$
\end{Thm}
 We have
$$ V(X,Y;L_1,L_2)=\hom(L_1,L_2)_{Y[1]\oplus X}$$ There is a natural $\bbc^{*}$-action on $Hom(L_1,L_2)_{Y[1]\oplus X}$
or $\mv(X,Y;L_1,L_2)$ simply by $t.(f,g,h)^{*}=(f,tg,h)^{*}$ for
$t\in \bbc^{*}$ and $(f,g,h)^{*}\in \mv(X,Y;L_1,L_2).$ We also
have a projective version of Theorem \ref{associativity}, where
$\mathbb{P}$ indicates the corresponding orbit space under the
$\bbc^*$-action.
\begin{Thm}\label{associativity2}
For fixed $A$-modules $X_i,Y_i,L_i$ for $i=1,2,$ we have
$$
\int_{\overline{Y}}g_{Y_2Y_1}^{Y}\chi(\mathbb{P}Hom(L_1,L_2)_{Y[1]\oplus
X})=\int_{\overline{L'}_1}g_{L'_1Y_1}^{L_1}\chi(\mathbb{P}Hom(L'_1,L_2)_{Y_2[1]\oplus
X}$$$$
\int_{\overline{X}}g_{X_2X_1}^{X}\chi(\mathbb{P}Hom(L_1,L_2)_{Y[1]\oplus
X})=\int_{\overline{L'}_2}g_{X_2L'_2}^{L_2}\chi(\mathbb{P}Hom(L_1,L'_2)_{Y[1]\oplus
X_1})
$$
\end{Thm}
\bigskip
\subsection{}We illustrate Theorem \ref{clustertheorem} by the following
example for infinite type quiver.
\begin{Example}
Let $Q$ be the Kronecker quiver $\xymatrix{1\ar @<2pt>[r]
\ar@<-2pt>[r]& 2}.$ Let $S_1$ and $S_2$ be simple modules of
vertex 1 and 2, respectively. Hence,
$$
R=\left(%
\begin{array}{cc}
  0 & 2 \\
  0 & 0 \\
\end{array}%
\right) \quad \mbox{and}\quad R'=\left(%
\begin{array}{cc}
  0 & 0 \\
  2 & 0 \\
\end{array}%
\right)
$$
$$
X_{S_1}=x^{\mathrm{\underline{dim}}S_1R'-\mathrm{\underline{dim}}S_1}+x^{\mathrm{\underline{dim}}S_1R-\mathrm{\underline{dim}}S_1}=x_1^{-1}(1+x_2^2)
$$
$$
X_{S_2}=x^{\mathrm{\underline{dim}}S_2R'-\mathrm{\underline{dim}}S_2}+x^{\mathrm{\underline{dim}}S_2R-\mathrm{\underline{dim}}S_2}=x_2^{-1}(1+x_1^2)
$$
For $\lambda\in \mathbb{P}^1(\bbc),$ let $u_{\lambda}$ be the
regular representation $\xymatrix{\bbc\ar @<2pt>[r]^{1}
\ar@<-2pt>[r]_{\lambda}& \bbc}.$ Then
$$
X_{u_{\lambda}}=x^{(1,1)R'-(1,1)}+x^{(1,1)R-(1,1)}+x^{(0,1)R+(1,0)R'-(1,1)}=x_1x_2^{-1}+x_1^{-1}x_2
+x_1^{-1}x_2^{-1}$$ Let $I_1$ and $I_2$ be indecomposable
injective modules corresponding vertex 1 and 2, respectively, then
$$
X_{(I_1\oplus I_2)[-1]}:=x^{\mathrm{\underline{dim}soc}(I_1\oplus
I_2)}=x_1x_2
$$
The left side of the identity of Theorem \ref{clustertheorem} is
$$
\mathrm{dim}_{\bbc}\mathrm{Ext}^{1}(S_1,S_2)X_{S_1}X_{S_2}=2(x_1^{-1}x_2^{-1}+x_1x_2^{-1}+x_1^{-1}x_2+x_1x_2)
$$
The first term of the right side is
$$
\int_{\lambda\in
\mathbb{P}^1(\bbc)}\chi(\mathbb{P}Ext^1(S_1,S_2)_{u_{\lambda}})X_{u_{\lambda}}=2(x_1^{-1}x_2^{-1}+x_1x_2^{-1}+x_1^{-1}x_2)
$$
As for the second term of the right side, we note that for any
$f\neq 0\in \mathrm{Hom}(S_2,\tau S_1),$ we have the exact
sequence:
$$
0\rightarrow S_2\xrightarrow{f}\tau S_1\rightarrow I_1\oplus
I_2\rightarrow 0
$$
This implies $\mathrm{Hom}(S_2,\tau S_1)\setminus
\{0\}=\mathrm{Hom}(S_2,\tau S_1)_{I_1\oplus I_2}$ Hence, the
second term is equal to $2x_1x_2.$
\end{Example}
\bigskip
\bigskip
\section{\bf 2-Calabi-Yau and Cluster multiplication formula}

A category $\mathcal{C}$ is 2-Calabi-Yau if there is an almost
canonical non degenerate bifunctorial pairing
$$\phi: \ext^1_{\mc}(M,N)\times \ext^{1}_{\mc}(N,M)\rightarrow \bbc
$$
for any object $M,N\in \mc.$ In particular, if $\mc$ is the
category of  nilpotent $\Lambda$-modules for some algebra
$\Lambda$ over $\bbc,$ then $\Lambda$ is called a 2-Calabi-Yau
algebra. In particular, if $\Lambda$ is the preprojective algebra
associated to a connected quiver $Q$ without loops, then $\Lambda$
is a 2-Calabi-Yau algebra. Let $\Lambda=\bbc Q/<R>$ be the
2-Calabi-Yau algebra associated to a connected quiver $Q$ without
loops and let $I$ be the set of vertices of $Q.$

Let $\Lambda_{\ud}$ be the affine variety of nilpotent
$\Lambda$-modules of dim.vector=$\ud.$ For
$\textbf{c}=(c_1,\cdots,c_m)\in \{0,1\}^m$ and $\textbf{j} =
(j_1,\cdots ,j_m)$ a sequence of elements of $I,$ let $x\in
\Lambda_{\ud}$, we define a $x$-stable flag of type
$(\textbf{j},\textbf{c})$ as a composition series of $x$
$$
\mathfrak{f}_{x} = \left(V=(\bbc^{\ud},x) \supseteq V^1 \supseteq
\cdots \supseteq V^m = 0\right)
$$
of $\Lambda$-submodules of $V$ such that $ |V^{k-1}/V^k|
=c_kS_{j_k} $ where $S_{j_k}$ is the simple module at the vertex
$j_k.$ Let $d_{\textbf{j},\textbf{c}}$ be the
 constructible functions satisfying
 $d_{\textbf{j},\textbf{c}}(x)=\chi(\Phi_{\textbf{j},\textbf{c},x})$ where $x\in \Lambda_{\ud}$ and
 $\Phi_{\textbf{j},\textbf{c},x}$ is the variety of $x$-stable flags of type $(\textbf{j},\textbf{c}).$
We simply write $d_\textbf{j}$ if $\textbf{c}=(1,\cdots,1).$
 Define $\mathcal{M}(\ud)$ to be the vector space spanned by
$d_{\textbf{j}}.$ Define
$$\Phi_{\textbf{j}}(\Lambda_{\ud})=\{(x,\mathfrak{f})\mid x\in
\Lambda_{\ud},\mathfrak{f}\in \Phi_{\textbf{j},x}\}.$$ We consider
a natural projection: $p:
\Phi_{\textbf{j}}(\Lambda_{\ud})\rightarrow \Lambda_{\ud},$ the
function $p_{*}(1_{\Phi_{\textbf{j}}(\Lambda_{\ud})})$ is
constructible by Theorem \ref{Joyce}.
\begin{Prop}
For a type $\textbf{j},$ the function $\Lambda_{\ud}\rightarrow
\bbq$ mapping $x$ to $\chi(\Phi_{\textbf{j},x})$ is constructible.
\end{Prop}
For fixed $\Lambda_{\ud},$ there are finitely many types
$\textbf{j}$ such that $\Phi_{\textbf{j}}(\Lambda_{\ud})$ is not
empty. Hence, there exists a finite subset $S(\ud)$ of
$\Lambda_{\ud}$ such that
$$\Lambda_{\ud}=\bigcup_{M\in S(\ud)}\str{M}$$ where $\str{M}=\{M'\in \Lambda_{\ud}\mid \chi(\Phi_{\textbf{j},M'})=\chi(\Phi_{\textbf{j},M}) \mbox{ for any type }\textbf{j}\}.$

For any $M\in \Lambda_{\ud},$ we define the evaluation form
$\delta_{M}:\mathcal{M}(\ud)\rightarrow \bbc$ mapping a
constructible function $f\in \mathcal{M}(\ud)$ to $f(M).$ We have
$$
\str{M}=\{M'\in \Lambda_{\ud}\mid \delta_{M'}=\delta_M\}
$$

We have the following multiplication formula \cite{XX2007b}.
\begin{Thm}\label{formula2}
For any $\Lambda$-modules $M$ and $N,$ we have
$$
\chi(\mathbb{P}\mathrm{Ext}^{1}_{\Lambda}(M,N))\delta_{M\oplus
N}=\sum_{L\in
S(\ue)}(\chi(\mathbb{P}\mathrm{Ext}^{1}_{\Lambda}(M,N)_{\str{L}})+\chi(\mathbb{P}Ext^{1}_{\Lambda}(N,M)_{\str{L}}))\delta_{L}
$$
where $\ue=\udim M+\udim N.$
\end{Thm}
\nd {\bf Remark.} In the case that $\Lambda$ is a preprojective
algebra, the above multiplication formula is given in
\cite[Theorem 1]{GLS2006}.

The proof of this theorem heavily depends on the fact that the
property of 2-Calabi-Yau induces the symmetry of the following two
linear maps. Define \cite [Section 2]{GLS2006}
$$
\beta_{\textbf{j},\textbf{c}',\textbf{c}'',\mathfrak{f}_{M},\mathfrak{f}_{N}}:
\bigoplus_{k=0}^{m-2}\ext^1_{\Lambda}(N_k,M_{k+1})\rightarrow
\bigoplus_{k=0}^{m-2}\ext^1_{\Lambda}(N_k,M_{k})
$$
by the following map
$$
\xymatrix{N_{k}\ar[r]^-{\varepsilon_k}\ar[d]^-{\iota_{N,k}}\ar[dr]&M_{k+1}[1]\ar[d]^-{\iota_{M,k+1}}\\N_{k-1}\ar[r]^-{\varepsilon_{k-1}}&M_{k}[1]}
$$
satisfying
$$
\beta_{\textbf{j},\textbf{c}',\textbf{c}'',\mathfrak{f}_{M},\mathfrak{f}_{N}}(\varepsilon_0,\cdots,\varepsilon_{m-2})=\iota_{M,1}\circ\varepsilon_0+\sum_{k=1}^{m-2}(\iota_{M,k+1}\circ\varepsilon_k-\varepsilon_{k-1}\circ
\iota_{N,k}).
$$
Depending on the 2-Calabi-Yau property of $\Lambda,$ we can write
down its dual.
$$
\beta'_{\textbf{j},\textbf{c}',\textbf{c}'',\mathfrak{f}_{M},\mathfrak{f}_{N}}:
\bigoplus_{k=0}^{m-2}\ext^1_{\Lambda}(M_k,N_{k})\rightarrow\bigoplus_{k=0}^{m-2}\ext^1_{\Lambda}(M_{k+1},N_{k})
$$
by the following map
$$
\xymatrix{M_{k+1}\ar[r]^-{\eta_{k+1}}\ar[d]^-{\iota_{M,k+1}}\ar[dr]&N_{k+1}[1]\ar[d]^-{\iota_{N,k+1}}\\M_{k}\ar[r]^-{\eta_k}&N_{k}[1]}
$$
satisfying
$$
\beta'_{\textbf{j},\textbf{c}',\textbf{c}'',\mathfrak{f}_{M},\mathfrak{f}_{N}}(\eta_0,\cdots,\eta_{m-2})=\sum_{k=0}^{m-3}(\eta_k\circ\iota_{M,k+1}-\iota_{N,k+1}\circ\eta_{k+1}
)+\eta_{m-2}\circ \iota_{M,m-1}.
$$
By the property of 2-Calabi-Yau, we have
$$
\mathrm{Ker}(\beta'_{\textbf{j},\textbf{c}',\textbf{c}'',\mathfrak{f}_{M},\mathfrak{f}_{N}})=\mathrm{Im}(\beta_{\textbf{j},\textbf{c}',\textbf{c}'',\mathfrak{f}_{M},\mathfrak{f}_{N}})^{\perp}.
$$
Of course, the next question is we should try to set up the same
formula over more general 2-Calabi-Yau categories.


\begin{thebibliography}{MRY90}

\bibitem[Bong]{Bong} K. Bongartz, \emph{Some geometric aspects of representation
theory}, In Algebras and modules, I (Trondheim, 1996), 1-27, CMS
Conf. Proc., 23, Amer. Math. Soc., Providence, RI, 1998. Solberg,
Eds.), Canad. Math. Soc. Conf. Proc. Series 23 (1998)

\bibitem[Borch]{Borch} R. Borcherds, \emph{Generalized Kac-Moody algebras}, J.
Algebra 115 (1988), 501-512.

\bibitem[BMRRT]{BMRRT} A. Buan, R. Marsh, M. Reineke, I. Reiten
and G. Todorov, \emph{Tilting theory and cluster combinatorics.}
Advances in Math. 204 (2006), 572-618.
\bibitem[CC]{CCformula}
P.~Caldero, F. Chapoton, \emph{Cluster algebras as Hall algebras
of quiver representations},  Commentarii Mathematici Helvetici 81
(2006), 595-616.
\bibitem[CCS]{CCS}
P.~Caldero, F. Chapoton, R. Schiffler, \emph{Quivers with
relations arising from clusters (A$_n$ case)}, Trans. Amer. Math.
Soc. 358 (2006),
 1347-1364.
\bibitem[CK1]{CK2005}
P. Caldero and B. Keller, \emph{From triangulated categories to
cluster algebras}, arXiv:math.RT/0506018. To appear in
Invent.math.

\bibitem[CK2]{CK2} P. Caldero and B. Keller, \emph{From triangulated categories to
cluster algebras II}, arXiv:math.RT/0510251v2, Annales Sci. del
l'Ecole Normale Super. 39(4): 83-100, 2006.

\bibitem[DX1]{DX1} B. Deng and J. Xiao, \emph{On double Ringel-Hall algebras}, J.
Algebra 251 (2002), 110-149.

\bibitem[DX2]{DX2} B. Deng and J. Xiao, \emph{On Ringel-Hall algebras}, Fields
Institute Communications Volume 40, 2004, 319-348.
\bibitem[Di]{Dimca}
 A. Dimca, \emph{Sheaves in topology.} Universitext. Springer-Verlag, Berlin, 2004.

\bibitem[DXX]{DXX}
M. Ding, J. Xiao and F. Xu, \emph{Realizing Enveloping Algebras via
Varieties of Modules}, math.QA/0604560.
\bibitem[FMV]{FMV} I. Frenkel, A. Malkin, and M. Vybornov, \emph{Affine Lie algebras and tame quivers}, Selecta Math.(N.S.) 2001, 7:1-56.
\bibitem[FZ1]{FZ2002}
S. Fomin, A. Zelevinsky, \emph{Cluster algebras. I. Foundations.}
J. Amer. Math. Soc. 15 (2002), no. 2, 497-529.

\bibitem[FZ2]{FZ2} S. Fomin, A. Zelevinsky, \emph{Cluster algebras. II. Finite type
classification.} Invent. Math. 154 (2003), no. 1, 63-121.

\bibitem[FZ3]{FZ3} S. Fomin, A. Zelevinsky, \emph{Cluster algebras: notes for the
CDM-03 conference}. Current developments in mathematics, 2003,
1-34, Int. Press, Somerville, MA, 2003.

\bibitem[Gre]{Green}
J.A. Green, \emph{Hall algebras, hereditary algebras and quantum
groups}, Invent. Math. 120, 361-377 (1995).

\bibitem[GLS]{GLS2006}C. Gei$\beta$, B. Leclerc, J. Schr\"oer, \emph{Semicanonical bases and
preprojective algebras II: A multiplication formula}, To appear in
Compositio Mathematica.

\bibitem[Hu1]{Hubery2004}
Andrew Hubery, \emph{From triangulated categories to Lie algebras:
A theorem of Peng and Xiao}, Proceedings of the Workshop on
Representation Theory of Algebras and related Topics (Quer¨¦taro,
2004), editors J. De la Pe\~na and R. Bautista.
\bibitem[Hu2]{Hubery2005}
Andrew Hubery, \emph{Acyclic cluster algebras via Ringel-Hall
algebras}, preprint.

\bibitem[Hu3]{Hubery2007}
Andrew Hubery, \emph{Hall Polynomials for Affine Quivers},
preprint.


\bibitem[Joy1]{Joyce-Hall} Joyce D, \emph{Configurations in abelian categories. II.
Ringel--Hall algebras},  Advances in Mathematics. 2007,
210:635-706.
\bibitem[Joy2]{Joyce} Joyce D, \emph{Constructible functions on Artin stacks}, J. London Math. Soc. 74, 583-606
(2006).

\bibitem[Ka]{Ka} M. Kapranov,\emph{Heisenberg doubles and derived categories}, J. Algebra 202 (1998), 712{ 744.
\bibitem[Ke1]{Keller1} B.~Keller, \emph{Chain complexes and stable categories}, Manus. Math. 67 (1990), 379--417.
\bibitem[Ke2]{keller2} B.~Keller, \emph{On triangulated orbit categories}.  Doc. Math. 10 (2005), 551-581.
\bibitem[KV]{KV} M. Kapranov and E. Vasserot, {\it Kleinian
singularities, derived categories and Hall algebras}, Math. Ann.
316 (2000), 565--576.
\bibitem[Lu1]{Lusztig} Lusztig G.  \emph{Quivers,perverse sheaves, and quantized enveloping algebras}, J. Amer. Math. Soc. 1991,
4(2):365-421.
\bibitem[Lu2]{Lusztig2} G. Lusztig, \emph{Semicanonical bases arising from enveloping algebras}, Adv. Math. 151 (2000),
no.2, 129-139.
\bibitem[Mac]{Macpherson} MacPherson R D, \emph{Chern classes for singular algebraic varieties}, Ann. Math. 100,
423-432 (1974).
\bibitem[Na]{Na} H. Nakajima, {\it Instantons on ALE spaces, quiver
varieties, and Kac--Moody algebras}, Duke Math. J. 76 (1994),
365--416.
\bibitem[Pe]{Peng} L. Peng, \emph{Lie algebras determined by Auslander-Reiten
quivers}, Comm. Algebra 26 (1998), no. 9, 2711-2725.
\bibitem[Rie]{Riedtmann} Ch. Riedtmann, \emph{Lie algebras generated by indecomposables}, J. Algebra 170,
526-546(1994).

\bibitem[Rin1]{Ringel} C. M. Ringel, \emph{ Hall algebras and quantum
groups,} Invent. Math. 101, 583-592 (1990).
\bibitem[Rin2]{RingelGreen}
C. M. Ringel, \emph{Green's theorem on Hall algebras},
Representation theory of algebras and related topics (Mexico City,
1994), 185--245, CMS Conf. Proc. 19, Amer. Math. Soc., Providence,
RI, (1996).
\bibitem[Ro]{Ro} M. Rosenlicht, \emph{A Remark on quotient spaces, spaces}, An. Acad.
Brasil. Ci$\hat{e}$nc. 35, 487-489 (1963).

\bibitem[Sc]{Schiffmann} O. Schiffmann, \emph{Lecture on Hall
algebras}, arXiv: math.RT/0611617.

\bibitem[SV]{SV} B. Sevenhant and M. Van den Bergh, \emph{A relation between a
conjecture of Kac and the structure of the Hall algebra}, J. Pure
and Appl. Algebra 160 (2001), 319-332.
\bibitem[To]{Toen2005} B.~To\"en, \emph{Derived Hall algebras},   Duke Math. J. 135 (2006), no. 3, 587-615.

\bibitem[X]{Xiao1997} J.Xiao, \emph{Drinfeld double and Ringel-Green theory of Hall
algebra}, J.Algebra 190 (1997), 100-144.

\bibitem[XX1]{XX2006}J. Xiao, and F. Xu, \emph{Hall algebras associated to triangulated categories},
preprint.
\bibitem[XX2]{XX2007a}J. Xiao, and F. Xu, \emph{Green's formula with $\bbc^{*}$-action and Caldero-Keller's formula for cluster
algebras}, preprint.
\bibitem[XX3]{XX2007b}J. Xiao, and F. Xu, \emph{A multiplication formula for 2-Calabi-Yau
algebras}, preprint.
\bibitem[XXZ]{XXZ2006}J. Xiao, F. Xu and G. Zhang, \emph{Derived categories and Lie algebras},
arXiv:math.QA/0604564.
\bibitem[ZW]{ZW}G. Zhang and S. Wang, \emph{The Green formula and heredity of algebras}, Science in
China, Ser. A Mathematics Vol 48, No.5 (2005), 610-617.}
\end{thebibliography}
\end{document}